\documentclass{amsart}

\usepackage{lmodern}
\usepackage[T1]{fontenc}
\usepackage[utf8]{inputenc}
\usepackage[english]{babel}
\usepackage{microtype} 

\usepackage{amsmath,amssymb,amsthm,mathrsfs,latexsym,mathtools,mathdots,enumerate,tikz}
\usetikzlibrary{arrows}
\usepackage[enableskew,vcentermath]{youngtab}
\usepackage{ytableau}

\usepackage[capitalize]{cleveref}

\newtheorem{theorem}{Theorem}
\newtheorem{proposition}[theorem]{Proposition}
\newtheorem{lemma}[theorem]{Lemma}
\newtheorem{corollary}[theorem]{Corollary}

\theoremstyle{definition}
\newtheorem{definition}[theorem]{Definition}
\newtheorem{example}[theorem]{Example}
\newtheorem{remark}[theorem]{Remark}
\newtheorem{question}[theorem]{Question}

\newcommand{\defin}[1]{\emph{#1}}

\newcommand{\setN}{\mathbb{N}}
\newcommand{\setZ}{\mathbb{Z}}
\newcommand{\setQ}{\mathbb{Q}}
\newcommand{\setR}{\mathbb{R}}

\newcommand{\Gset}{\mathbf{G}}
\newcommand{\polyP}{\mathcal{P}}

\newcommand{\xvec}{\mathbf{x}}
\newcommand{\yvec}{\mathbf{y}}

\newcommand{\wvec}{\mathbf{w}}

\newcommand{\onevec}{\mathbf{1}}
\newcommand{\lambdavec}{{\boldsymbol\lambda}}
\newcommand{\tauvec}{{\boldsymbol\tau}}
\newcommand{\muvec}{{\boldsymbol\mu}}
\newcommand{\nuvec}{{\boldsymbol\nu}}

\newcommand{\tinsert}{\boxtimes} 

\newcommand{\tiling}{\mathscr{P}}

\newcommand{\refin}{\leq_{\mathrm{ref}}}
\newcommand{\refiner}{<_{\mathrm{ref}}}
\newcommand{\interl}{\geq_{\mathrm{int}}}

\DeclareMathOperator*{\length}{\mathit{l}}

\title[Gelfand--Tsetlin polytopes]{Gelfand--Tsetlin polytopes and the integer decomposition property}

\author[P.~Alexandersson]{Per Alexandersson}
\email{per.w.alexandersson@gmail.com}

\begin{document}

\begin{abstract}

Let $\polyP$ be the Gelfand--Tsetlin polytope
defined by the skew shape $\lambdavec/\muvec$ and weight $\wvec$.
In the case corresponding to a standard Young tableau,
we completely characterize for which shapes $\lambdavec/\muvec$
the polytope $\polyP$ is integral.
Furthermore, we show that $\polyP$ is a compressed polytope
whenever it is integral and corresponds to a standard Young tableau.
We conjecture that a similar property hold for
arbitrary $\wvec$, namely that $\polyP$ has the integer decomposition property
whenever it is integral. 

Finally, a natural partial ordering on GT-polytopes is introduced
that provides information about integrality and the integer decomposition property,
which implies the conjecture for certain shapes.
\medskip

\end{abstract}

\maketitle

\setcounter{tocdepth}{1}
\tableofcontents

\section{Introduction}

The study of polytopes related to quantities in representation
theory has been fruitful over the last few decades; a highlight, for example,
is Knutson and Tao's proof of the Saturation conjecture \cite{Knutson99thehoneycomb}.

There are two types of (skew) Gelfand--Tsetlin polytopes, which are defined later on.
The first type is weight-restricted polytopes $\polyP_{\lambdavec/\muvec,\wvec}$,
which is the main topic of this paper.
There are also polytopes without a restriction on the weight,
$\polyP_{\lambdavec/\muvec}$, which are only covered briefly.

In order to study linear recurrence relations among skew Schur polynomials,
the integer decomposition property of Gelfand--Tsetlin polytopes
without weight restriction was proved and used in \cite{Alexandersson20141}.
A sketch of this proof can be found in \cref{sec:noweight}.

It is therefore natural to consider the weight-restricted version of Gelfand--Tsetlin polytopes.
These polytopes are more complicated;
for example, not all such polytopes are integral, as proved by De Loera and McAllister in \cite{Loera04}.
\medskip 

A skew Gelfand--Tsetlin polytope $\polyP_{\lambdavec/\muvec,\wvec}$ is
defined by two partitions, $\lambdavec/\muvec$ and an integer composition $\wvec$
(the exact definition is given in \cref{sec:gtpolytopes}).

\smallskip

\noindent 
\textbf{Main result:}
The main result of this paper (\cref{onepartitioniscompressed}) concerns the case
$\wvec=(1,\dots,1)$, which correspond to the case of \emph{standard} Young tableaux.
We completely characterize the skew shapes $\lambdavec/\muvec$ for which
 $\polyP_{\lambdavec/\muvec,\wvec}$ is integral, see \cref{cor:completeIntegralitycharacterization},
and then proceed to show that each such integral polytope is \emph{compressed}.

\smallskip

\noindent 
\textbf{Second result:}
We show that if $\wvec'$ is a refinement of $\wvec$, then 
\begin{itemize}
\item If $\polyP_{\lambdavec/\muvec,\wvec'}$ is integral, then $\polyP_{\lambdavec/\muvec,\wvec}$ is integral.
\item If $\polyP_{\lambdavec/\muvec,\wvec'}$ has the integer decomposition property,
then $\polyP_{\lambdavec/\muvec,\wvec}$ is also has the integer decomposition property.
\end{itemize}
This result is presented as \cref{thm:paoprop}
and it implies integrality and non-integrality of several natural
families of Gelfand--Tsetlin polytopes.
For example, all hook shapes and disjoint unions of rows yield integral polytopes,
which also have the integer decomposition property.
The latter family of polytopes has a natural interpretation in terms of
contingency matrices, (see \cref{sec:contingencyMatrices}).
We conjecture that $\polyP_{\lambdavec/\muvec,\wvec}$ is integral if and only if it
has the integer decomposition property. This is supported by computer experiments.

The present article extends previous work by
King, Tollu, Toumazet \cite{King04stretched}, De Loera and McAllister \cite{Loera04,LoeraM06};
we extensively use and extend results by the latter two authors.

The work \cite{RassartThesis} by Rassart gives a good overview on
the connection between representation theory and polytopes.
Briefly stated, there is a bijection (given further down) between integral points 
inside Gelfand--Tsetlin polytopes and semi-standard
Young tableaux, which explains the connection with representation theory.

\subsection*{Acknowledgements}

The author would like to thank Valentin Féray, Christian Haase, Tyrrell McAllister and Benjamin Nill 
for helpful discussions, as well as the anonymous referees for their suggestions.
This work has been partially funded by the Knut and Alice Wallenberg Foundation.

\section{Preliminaries}

We expect that the reader is somewhat familiar with the notion of Young tableaux and skew Young tableaux.
A standard reference in this field is \cite{Macdonald79symmetric}.
\medskip

Let $\lambdavec$ and $\muvec$ be integer partitions where $\lambda_i \geq \mu_i$.
A \defin{skew Young diagram} of \defin{shape} $\lambdavec/\muvec$ is an 
arrangement of ``boxes'' in the plane with coordinates given by
$\{ (i,j) \in \setZ^2 | \mu_i < j \leq \lambda_i \}$. Note that we use the English convention.
For example, the skew diagram of shape $\lambdavec/\muvec=(5,4,2,2)/(2,1)$ is given to the left in \eqref{eq:diagramandtableau}.
\begin{equation}\label{eq:diagramandtableau}
\young(::\hfil\hfil\hfil,:\hfil\hfil\hfil,\hfil\hfil,\hfil\hfil) \qquad
\young(::112,:123,13,44).
\end{equation}
A \defin{semi-standard Young tableau} (or SSYT) is a Young diagram with natural numbers in the boxes,
such that each row is weakly increasing and each column is strictly increasing, as in \eqref{eq:diagramandtableau}.
Whenever all numbers are different, we say that the tableau is \defin{standard},
and whenever $\muvec$ is not the empty partition, it is a \defin{skew} tableau.
\medskip

\subsection{Notation}

We always use bold lowercase letters, $\xvec$, to denote vectors or partitions $(x_1,x_2,\dots,x_n)$.
The symbol $\onevec$ denotes the integer composition $(1,1,\dots,1)$
where the length is evident from the context. The sum of the entries in the vector $\xvec$ is denoted $|\xvec|$.
If $S$ is a set, $|S|$ also denotes the cardinality of $S$.
Sets are always capital letters, so there should be no confusion.
Whenever we need multiple vectors, we index these with superscript. Multiplication by a constant
is done elementwise on vectors, integer compositions, and partitions.
The number $\length(\wvec)$ denote the index of the last non-zero entry in $\wvec$.

We write $\lambdavec \geq \muvec$ if $\lambda_i \geq \mu_i$ for all $i$.
A stronger condition is defined by the partial order $\lambdavec \interl \muvec$,
which indicates that the entries in $\lambdavec$ and $\muvec$ \defin{interlace}:
\[
\lambda_1\geq \mu_1 \geq \lambda_2 \geq \mu_2 \geq \dots \geq \lambda_n \geq \mu_n.
\]
Clearly, $\lambdavec \interl \muvec$ if and only if $k\lambdavec \interl k\muvec$ for $k>0$.
For integer compositions $\wvec$ and $\wvec'$, we write $\wvec' \refin \wvec$ 
if $\wvec'$ is a composition refinement of $\wvec$.

As per the standard convention, partitions are padded with zeros as needed.

\subsection{GT-patterns}

Gelfand--Tsetlin patterns were first introduced in \cite{GelfandTsetlin50}.
A \defin{Gelfand--Tsetlin pattern}, or GT-pattern for short, is a
triangular or parallelogram arrangement of non-negative numbers,
\[
\begin{array}{cccccccccccccc}
x^m_{1} & & x^m_{2} & & \cdots & & \cdots & & x^m_{n} \\
 & \ddots & & \ddots &  & &   & & & \ddots  \\
  &  &   x^2_{1} &  & x^2_{2} & & \cdots & & \cdots & & x^2_{n} \\
   &  &    & x^1_{1} &    & x^1_{2} & & \cdots & & \cdots & & x^1_{n}
\end{array}
\]
where the entries must satisfy the inequalities
\begin{equation}
x^{i+1}_{j} \geq x^i_{j} \text{ and } x^i_{j} \geq x^{i+1}_{j+1} \label{eq:gtinequalities}
\end{equation}
for all values of $i$, $j$ where the indexing is defined.
The inequalities simply states that horizontal rows and down-right diagonals are weakly decreasing,
while down-left diagonals are weakly increasing.
Note that the conditions ensure that any two adjacent rows interlace; $\xvec^{i+1} \interl \xvec^{i}$.
\medskip 

Whenever all $x^i_j$ are natural numbers, we say that the GT-pattern is \defin{integral}.
There is a bijection between integral GT-patterns and skew Young tableaux;
The skew shape defined by row $j$ and $j+1$ in a GT-pattern $G$
describes which boxes in a tableau $T$ that have content $j$.
In particular, if the bottom row in $G$ is $\muvec$ and the top row is $\lambdavec$,
then $T$ has shape $\lambdavec/\muvec$.
Here is an example of this correspondence:
\begin{equation}\label{eq:gtcorr}
\setcounter{MaxMatrixCols}{20}
\begin{matrix}
4 &   & 3 &   & 2 &   & 1\\
 & 4 &   & 3 &   & 1 &   & 1\\
 &  & 3 &   & 3 &   & 1 &   & 1\\
 &  &  & 3 &   & 2 &   & 1 &   & 0\\
 &  &  &  & 2 &   & 1 &   & 1 &   & 0\\
\end{matrix}
\quad
\longleftrightarrow
\quad
\young(::13,:12,:4,2)
\end{equation}

\medskip

Note that if $|\muvec|=0$, then $x^i_j=0$ whenever $j \geq i$.
In this case, these entries are usually not displayed and it suffices to present a triangular array,
which is the more common form of GT-patterns;
\begin{equation}\label{eq:triangulargt}
\begin{matrix}
5 &   & 4 &   & 2 &   & 1 &   & 1 &  & 0\\
 & 5 &   & 3 &   & 2 &   & 1 &   & 0\\
 &  & 3 &   & 3 &   & 2 &   & 1  \\
 &  &  & 3 &   & 3 &   & 1 &  \\
 &  &  &  & 3 &   & 2 &  \\
 &  &  &  &  & 3 & \\
\end{matrix}
\quad
\longleftrightarrow
\quad
\young(11155,2236,34,4,6)
\end{equation}
It is customary to only write $\lambdavec$ instead of $\lambdavec/\muvec$ to describe such a shape.

\begin{remark}\label{rem:gtboxcount}
In any GT-pattern, $x^{i+1}_{j} - x^i_{j}$ counts the
number of boxes with content $i$ in row $j$ in the corresponding tableau.
\end{remark}

Let $\Gset_{\lambdavec/\muvec,\wvec}$ denote the set of all integral GT-patterns
where the top row is given by $\lambdavec$, the bottom row is given by $\muvec$ and
\begin{equation}
\sum_{j=1}^n (x^{i+1}_{j} - x^i_{j}) = |\xvec^{i+1}|-|\xvec^{i}| = w_i \text{ for } i=1,2,\dotsc,\length(\wvec). \label{eq:gtweightequalities}
\end{equation}
The integer composition $\wvec$ is usually referred to as the \defin{weight}.
We let $\Gset^k_{\lambdavec/\muvec,\wvec}$ be a short form of $\Gset_{k\lambdavec/k\muvec,k\wvec}$.
In most applications, $\wvec$ is always considered to be a partition,
since $|\Gset_{\lambdavec/\muvec,\wvec}|=|\Gset_{\lambdavec/\muvec,\pi(\wvec)}|$
for any permutation $\pi$. This identity not obvious, but follows from \cite{BenderKnuth1972}.

\subsection{Adding GT-patterns}\label{ssec:adding}

Given two GT-patterns $G_1 \in \Gset^{k_1}_{\lambdavec/\muvec,\wvec}$  and $G_2 \in  \Gset^{k_2}_{\lambdavec/\muvec,\wvec}$,
we define $G_1+G_2$ as element-wise addition on the patterns and it is easy to verify that
\begin{equation}\label{eq:gtpplus}
G_1 + G_2 \in  \Gset^{k_1 + k_2}_{\lambdavec/\muvec,\wvec}.
\end{equation}
As an example,
\[
\begin{matrix}
4 &   & 3 &   & 2\\
 & 4 &   & 3 &   & 1\\
 &  & 3 &   & 3 &   & 1\\
 &  &  & 3 &   & 2 &   & 0\\
 &  &  &  & 2 &   & 1 &   & 0\\
\end{matrix}
\mkern-30mu
+
\mkern-30mu
\begin{matrix}
8 &   & 6 &   & 4\\
 & 7 &   & 6 &   & 3\\
 &  & 7 &   & 4 &   & 3\\
 &  &  & 7 &   & 3 &   & 0\\
 &  &  &  & 4 &   & 2 &   & 0\\
\end{matrix}
\mkern-30mu
=
\mkern-30mu
\begin{matrix}
12 &   & 9 &   & 6\\
 & 11 &   & 9 &   & 4\\
 &  & 10 &   & 7 &   & 4\\
 &  &  & 10 &   & 5 &   & 0\\
 &  &  &  & 6 &   & 3 &   & 0\\
\end{matrix}.
\]
This operation (which we denote $\tinsert$) can also be seen on the Young tableaux side,
\[
\young(::13,:12,24)
\quad
\tinsert
\quad
\young(::::1114,::1233,2224)
\quad
=
\quad
\young(::::::111134,:::112233,222244).
\]
One can deduce from \eqref{eq:gtweightequalities} that $\tinsert$ corresponds
to row-wise concatenation of the corresponding
Young tableaux and rearrangement of the boxes in increasing order.

\medskip

It should be noted that there are few sources that cover GT-patterns of skew type. 
There are several reasons for this: there are not as many applications of skew Young tableaux in representation theory,
and many cases can be reduced to a non-skew setting.

However, it will be evident in the techniques used in this paper that it is less painful to work
with parallelogram rather than triangular arrangements.

\section{Gelfand--Tsetlin polytopes}\label{sec:gtpolytopes}

Given a skew shape $\lambdavec/\muvec$ and weight $\wvec$, the \defin{Gelfand--Tsetlin polytope} $\polyP_{\lambdavec/\muvec,\wvec} \subset \setR^{mn}$ (or GT-polytope),
is defined as the convex polytope consisting of all GT-patterns $(x^i_j)_{1\leq i \leq m,1\leq j \leq n}$ that satisfy 
the equalities 
\begin{itemize}
\item $\xvec^m = \lambdavec$ and  $\xvec^1 = \muvec$;
\item $|\xvec^{i+1}|-|\xvec^{i}| = w_i$ for $i=1,2,\dots,m-1$.
\end{itemize}
The elements in $\Gset_{\lambdavec/\muvec,\wvec}$ are in bijection with the
integer lattice points in $\polyP_{\lambdavec/\muvec,\wvec}$.
GT-polytopes (mainly those for which $|\muvec|=0$ and $\wvec$ is a partition) have been studied in
many places before, \emph{e.g.}~\cite{Kirillov88thebethe,Loera04,Rassart2004}.

In order to have consistent notation, let $\polyP^k_{\lambdavec/\muvec,\wvec}$ denote the $k$-dilation of the polytope $\polyP_{\lambdavec/\muvec,\wvec}$.
Analogous to \eqref{eq:gtpplus}, it is easy to verify the following identity on the Minkowski sum of GT-polytopes:
\[
\polyP^{k_1}_{\lambdavec/\muvec,\wvec} + \polyP^{k_2}_{\lambdavec/\muvec,\wvec} = \polyP^{k_1 + k_2}_{\lambdavec/\muvec,\wvec}.
\]

A convex polytope is called \defin{integral} if all of its vertices are integer points.
Gelfand--Tsetlin polytopes are in general \emph{not} integral, see \cite{King04stretched,Loera04},
but it is known that the Ehrhart quasi-polynomial for such a polytope is actually a polynomial;
this is not to be expected of a general non-integral polytope.
See \cite{Kirillov88thebethe,Rassart2004} for a proof of the polynomiality of the Ehrhart quasi-polynomial.

Even though $\Gset_{\lambdavec/\muvec,\wvec}$ and $\Gset_{\lambdavec/\muvec,\pi(\wvec)}$ always
have the same cardinality for any permutation $\pi$, the corresponding polytopes may look different.
For example, the parameters $\lambdavec=(5,3)$, $\muvec=(0)$ and $\wvec=(2,2,2,1,1)$ give an integral
polytope, but setting $\wvec=(2,2,1,2,1)$ for the same shape yields a non-integral polytope.
This explains the more general setting of composition weights.

\subsection{The geometry of GT-polytopes}

In this subsection, we recall some notions introduced in \cite{Loera04}.
These definitions and results were originally proved only for non-skew shapes,
but the same proofs can be carried out in the skew setting, which is what is stated here.

\begin{definition}[Tiling]
The \defin{tiling} of a GT-pattern $G$ is a partition $\tiling$ of entries in $G$ into \defin{tiles}.
This partition is defined as the finest partition with the property that
entries in $G$ that are equal and adjacent\footnote{Adjacent in the four directions NW, NE, SW and SE.} belong to the same tile.

Thus, $G$ is constant on each of its tiles and for each tile, this constant is called the \defin{content} of the tile.
Tiles that contain entries from the bottom or the top row are called \defin{fixed}, 
while all other tiles are \defin{free}.
\end{definition}

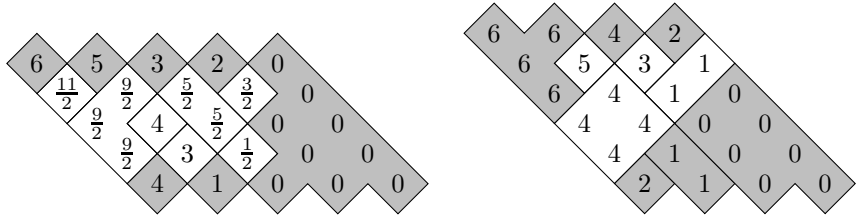
\begin{figure}[!ht]
\begin{tikzpicture}[scale=0.4]
\draw[black] (3,-1)--(4,-2)--(5,-1)--(4,0)--(3,-1)--cycle;
\draw[black] (6,-2)--(7,-1)--(6,0)--(5,-1)--(4,-2)--(5,-3)--(6,-4)--(7,-3)--(6,-2)--cycle;
\draw[black] (9,-1)--(8,0)--(7,-1)--(8,-2)--(9,-3)--(10,-2)--(9,-1)--cycle;
\draw[black] (9,-1)--(10,-2)--(11,-1)--(10,0)--(9,-1)--cycle;
\draw[black] (6,-2)--(7,-3)--(8,-2)--(7,-1)--(6,-2)--cycle;
\draw[black] (7,-3)--(8,-4)--(9,-3)--(8,-2)--(7,-3)--cycle;
\draw[black] (9,-3)--(10,-4)--(11,-3)--(10,-2)--(9,-3)--cycle;
\filldraw[color=black,fill=lightgray] (2,0)--(3,-1)--(4,0)--(3,1)--(2,0)--cycle;
\filldraw[color=black,fill=lightgray] (4,0)--(5,-1)--(6,0)--(5,1)--(4,0)--cycle;
\filldraw[color=black,fill=lightgray] (6,0)--(7,-1)--(8,0)--(7,1)--(6,0)--cycle;
\filldraw[color=black,fill=lightgray] (8,0)--(9,-1)--(10,0)--(9,1)--(8,0)--cycle;
\filldraw[color=black,fill=lightgray] (6,-4)--(7,-5)--(8,-4)--(7,-3)--(6,-4)--cycle;
\filldraw[color=black,fill=lightgray] (8,-4)--(9,-5)--(10,-4)--(9,-3)--(8,-4)--cycle;
\filldraw[color=black,fill=lightgray] (10,0)--(11,-1)--(10,-2)--(11,-3)--(10,-4)--(11,-5)--(12,-4)--(13,-5)--(14,-4)--(15,-5)--(16,-4)--(15,-3)--(14,-2)--(13,-1)--(12,0)--(11,1)--(10,0)--cycle;
\node at (3,0) {$6$};
\node at (5,0) {$5$};
\node at (7,0) {$3$};
\node at (9,0) {$2$};
\node at (11,0) {$0$};
\node at (4,-1) {$\frac{11}{2}$};
\node at (6,-1) {$\frac{9}{2}$};
\node at (8,-1) {$\frac{5}{2}$};
\node at (10,-1) {$\frac{3}{2}$};
\node at (12,-1) {$0$};
\node at (5,-2) {$\frac{9}{2}$};
\node at (7,-2) {$4$};
\node at (9,-2) {$\frac{5}{2}$};
\node at (11,-2) {$0$};
\node at (13,-2) {$0$};
\node at (6,-3) {$\frac{9}{2}$};
\node at (8,-3) {$3$};
\node at (10,-3) {$\frac{1}{2}$};
\node at (12,-3) {$0$};
\node at (14,-3) {$0$};
\node at (7,-4) {$4$};
\node at (9,-4) {$1$};
\node at (11,-4) {$0$};
\node at (13,-4) {$0$};
\node at (15,-4) {$0$};
\end{tikzpicture}
\quad
\begin{tikzpicture}[scale=0.4]
\draw[black] (5,-1)--(6,-2)--(7,-1)--(6,0)--(5,-1)--cycle;
\draw[black] (7,-1)--(8,-2)--(9,-1)--(8,0)--(7,-1)--cycle;
\draw[black] (8,-2)--(9,-3)--(10,-2)--(11,-1)--(10,0)--(9,-1)--(8,-2)--cycle;
\draw[black] (8,-4)--(9,-3)--(8,-2)--(7,-1)--(6,-2)--(5,-3)--(6,-4)--(7,-5)--(8,-4)--cycle;
\filldraw[color=black,fill=lightgray] (6,0)--(7,-1)--(8,0)--(7,1)--(6,0)--cycle;
\filldraw[color=black,fill=lightgray] (8,0)--(9,-1)--(10,0)--(9,1)--(8,0)--cycle;
\filldraw[color=black,fill=lightgray] (7,-5)--(8,-6)--(9,-5)--(8,-4)--(7,-5)--cycle;
\filldraw[color=black,fill=lightgray] (10,-4)--(9,-3)--(8,-4)--(9,-5)--(10,-6)--(11,-5)--(10,-4)--cycle;
\filldraw[color=black,fill=lightgray] (6,-2)--(5,-1)--(6,0)--(5,1)--(4,0)--(3,1)--(2,0)--(3,-1)--(4,-2)--(5,-3)--(6,-2)--cycle;
\filldraw[color=black,fill=lightgray] (10,-2)--(9,-3)--(10,-4)--(11,-5)--(12,-6)--(13,-5)--(14,-6)--(15,-5)--(14,-4)--(13,-3)--(12,-2)--(11,-1)--(10,-2)--cycle;
\node at (3,0) {$6$};
\node at (5,0) {$6$};
\node at (7,0) {$4$};
\node at (9,0) {$2$};
\node at (4,-1) {$6$};
\node at (6,-1) {$5$};
\node at (8,-1) {$3$};
\node at (10,-1) {$1$};
\node at (5,-2) {$6$};
\node at (7,-2) {$4$};
\node at (9,-2) {$1$};
\node at (11,-2) {$0$};
\node at (6,-3) {$4$};
\node at (8,-3) {$4$};
\node at (10,-3) {$0$};
\node at (12,-3) {$0$};
\node at (7,-4) {$4$};
\node at (9,-4) {$1$};
\node at (11,-4) {$0$};
\node at (13,-4) {$0$};
\node at (8,-5) {$2$};
\node at (10,-5) {$1$};
\node at (12,-5) {$0$};
\node at (14,-5) {$0$};
\end{tikzpicture}
\caption{The tilings of a non-integral and an integral GT-pattern. The gray tiles are the fixed tiles.}\label{fig:tilingexample}
\end{figure}

Two Gelfand--Tsetlin patterns are displayed with their tilings in Figure \ref{fig:tilingexample}.

\medskip

The following definition is slightly different, but for our purposes, equivalent with the definition in \cite{Loera04}:
\begin{definition}[Tiling matrix]
Let $G$ be a GT-pattern with $m$ rows and $s$ free tiles, enumerated in some way.
A \defin{tiling matrix} $T_G = (t_{ij})_{1 \leq i \leq m, 1\leq j \leq s}$ associated to $G$ is a matrix with $s$ columns and $m$ rows, 
such that $t_{ij}$ is the number of entries in the \emph{free tile} $j$ that are in the $(m-i+1)$th row  of $G$.
\end{definition}
Tiling matrices for the GT-patterns in Figure \ref{fig:tilingexample} are for example
\[
\begin{pmatrix}
 0 & 0 & 0 & 0 & 0 & 0 & 0 \\
 1 & 1 & 1 & 1 & 0 & 0 & 0 \\
 0 & 1 & 1 & 0 & 1 & 0 & 0 \\
 0 & 1 & 0 & 0 & 0 & 1 & 1 \\
 0 & 0 & 0 & 0 & 0 & 0 & 0 \\
\end{pmatrix}
\qquad 
\text{ and }
\qquad 
\begin{pmatrix}
 0 & 0 & 0 & 0 \\
 1 & 1 & 1 & 0 \\
 0 & 0 & 1 & 1 \\
 0 & 0 & 0 & 2 \\
 0 & 0 & 0 & 1 \\
 0 & 0 & 0 & 0 \\
\end{pmatrix}.
\]
In the first matrix, the columns corresponds to the free tiles with 
content $\frac{11}{2}$, $\frac{9}{2}$, $\frac{5}{2}$, $\frac{3}{2}$, $4$, $3$ and  $\frac12$.
The columns in the second matrix corresponds to the contents $5$, $3$, $1$ and $4$.
Note that the first and last row in any tiling matrix only contain zeros.
\medskip 

We can now state the main theorem in \cite{Loera04}:
\begin{theorem}\label{thm:mcallistermain}
Suppose $T_G$ is the tiling matrix of some $G \in \polyP_{\lambdavec/\muvec,\wvec}$.
Then the dimension of $\ker T_G$ is equal to the dimension of 
the minimal (dimensional) face of the GT-polytope containing $G$.
\end{theorem}
Note that $\dim \ker T_G $ 
is independent of the order of the columns in $T_G$.
Also note that a GT-pattern is a vertex of $\polyP_{\lambdavec/\muvec,\wvec}$ 
if and only if the columns in the tiling matrix are linearly independent.
\medskip

Let $\tiling$ and $\tiling'$ be tilings of GT-patterns in some $\polyP_{\lambdavec/\muvec,\wvec}$.
A tiling $\tiling'$ is a \defin{refinement} of $\tiling$
if each tile of $\tiling'$ is a subtile of a tile in $\tiling$.

\begin{lemma}\label{lem:tilingrefinement}
Let $G = a_1G_1 + a_2G_2 + \dots + a_kG_k$ with $G_j \in \polyP_{\lambdavec/\muvec,\wvec}$ and all $a_j>0$.
Then the tiling of $G$ is a refinement (not necessarily strict) of the tiling of each $G_j$.
\end{lemma}
\begin{proof}
All GT-patterns have the same monotonicity of the entries along down-left, down-right, up-left and up-right
diagonals, so if two adjacent entries in $G$ are equal, the corresponding entries in $G_j$ must be equal.
\end{proof}

The following lemma shows how the kernel of the tiling matrix relate a GT-pattern $G$ with
vertices of the minimal-dimensional face that contains $G$.
\begin{lemma} \label{lem:tilingkernel}
Let $G$ be a non-vertex of a GT-polytope $\polyP_{\lambdavec/\muvec,\wvec}$
and let $G'$ be a vertex of the minimal-dimensional face that contains $G$.
If $T_i$ is a tile of $G$, let $x_i$ be the content of tile $i$ in $G$ and let $y_i$ be the content of $T_i$ when viewed as a subset of $G'$.
\emph{The latter is well-defined, since according to \cref{lem:tilingrefinement}, $G'$ is constant on each tile of $G$.}

Then $\xvec-\yvec$ is in the kernel of $T_G$.
\end{lemma}
\begin{proof}
For each entry in $G$ that is a member of a fixed tile, the corresponding entry in $G'$ is also a member of
a fixed tile. Thus, $G$ and $G'$ agrees on the fixed tiles of $G$.
Let $\wvec'_i$ be the sum of the entries in the fixed tiles of $G$ in row $i$.
Adding all entries in each row of $G$ (free plus fixed) we have that $T_G \xvec + \wvec' = \wvec$.
Similarly, we have that $T_G \yvec + \wvec' = \wvec$ by adding all entries in each row of $G'$.
From this, it is evident that $T_G(\xvec-\yvec)=0$.
\end{proof}

\bigskip

We end this section with a small application of tiling matrices:
\begin{proposition}
If $\lambdavec$ is a hook, \emph{i.e.} of the form $\lambdavec = (h,1,1,\dots,1)$,
then for any $\wvec$, \emph{all} integral points in $\polyP_{\lambdavec,\wvec}$ are vertices of $\polyP_{\lambdavec,\wvec}$.
\end{proposition}
\begin{proof}
All integral GT-patterns of hook shape are of the following form:
\[
\begin{matrix}
h &&  1 && \dots  && 1  && 0  && \dots && 0\\
& \ddots && \ddots && \iddots && \star && \ddots && \iddots \\
 && x_j && 1 && \star && \star && 0 \\
 &&& \ddots && \ddots && \star && \iddots \\
 &&&& \ddots && \ddots && \iddots \\
 &&&&& x_2 && \star \\
 &&&&&& x_1
\end{matrix}
\]
Since NW-SE diagonals are decreasing, all $\star$ entries must be either $0$ or $1$.
These entries must therefore belong to the fixed tiles and the only free tiles consist of (possibly a subset of) the $x_i$.
Since no two such entries appear in the same row, the corresponding tiling
matrix have full rank and it follows that $G$ is a vertex.
\end{proof}

\subsection{Geometric properties of Gelfand--Tsetlin polytopes}

We now prove some results regarding integrality and non-integrality of some natural families of GT-polytopes.

\begin{lemma}\label{lem:unimodulartilingmatrix}
Let $G$ be a vertex of $\polyP_{\lambdavec/\muvec,\wvec}$ and suppose that the tiling matrix $T_G$ has rank $s$.
If every $s \times s$-minor of $T_G$ has determinant $\pm 1$ or $0$, then $G$ is integral.
\end{lemma}
\begin{proof}
Since $G$ is a vertex, the $s$ columns of $T_G$ are linearly independent.
Let $\xvec = (x_1,x_2,\dots,x_s)$ be the vector such that $x_j$ is the content of tile $j$.

The non-free tiles only contain integers and each row sum in $G$ is an integer.
The vector $\yvec = T_G \xvec$ must therefore be an integral vector;
$y_j$ is the sum of the entries in row $j$ in $G$ that belong to free tiles.

Since every $s \times s$-minor of $T_G$ has determinant $\pm 1$ or $0$ and at least one such minor is non-zero,
it follows that there is an invertible integer matrix $U$ such that $U \yvec = UT_G \xvec$
and the top $s \times s$-submatrix of $UT_G$ is the identity matrix.
This implies that entries in $\xvec$ are integers and hence $G$ is integral.
\end{proof}

\medskip

\begin{proposition}\label{prop:onevertices}
All integral GT-patterns in $\polyP_{\lambdavec/\muvec,\onevec}$
are vertices of $\polyP_{\lambdavec/\muvec,\onevec}$.
\end{proposition}
\begin{proof}
Consider an integral GT-pattern $G \in \polyP_{\lambdavec/\muvec,\onevec}$.

Let $T$ be a tiling matrix of $G$ and let $T_i$ be the submatrix of $T$,
given by the $i$ bottom rows of $T$ and the columns corresponding to
the free tiles that intersect at least one of the rows $1,2,\dots,i$ in $G$.
By definition, if $G$ has $m$ rows, then $T_m=T$.
\medskip

The goal is to show that the columns of $T_i$ are linearly independent for all $i$.
The case $T_1$ is clear, since this matrix has no columns.
Now assume that all columns of $T_i$ are linearly independent
and consider the rows $i$ and $i+1$ in $G$.
Since the row sum only increase by one from row $i$ to $i+1$
and $G$ is integral, there is only one entry which differs between adjacent rows.
Hence, if $\tauvec^i$ is row $i$, row $i$ and $i+1$ must have the form
\begin{equation*}
\begin{matrix}
\tau^i_1 &&  \dots && \tau^i_{j-1}   && \boxed{ \tau^i_{j}+1 }  &&\tau^i_{j +1}  &&\dots && \tau^i_n \\
&\tau^i_1 &&  \dots && \tau^i_{j-1}   &&\tau^i_{j}   &&\tau^i_{j +1}  &&\dots && \tau^i_n
\end{matrix}.
\end{equation*}
There are two cases to consider, depending on the value of the boxed entry:
\begin{align*}
\begin{cases}
\tau^i_{j}+1 = \tau^i_{j-1} & \Longrightarrow \quad  T_{i+1} = \begin{pmatrix} \star \\ T_i \end{pmatrix} \\
\tau^i_{j}+1 < \tau^i_{j-1} & \Longrightarrow \quad
T_{i+1} = \begin{pmatrix}
\star & 1 \\
T_i & 0
\end{pmatrix}
\end{cases}.
\end{align*}
\emph{Here, the entries in $\star$ is unimportant.  It is safe to assume that the matrices have this form,
as any permutation of the columns do not change the rank of the matrix.}

In the first case, the boxed entry is a member of a free tile which intersects row $i$,
so it does not contribute to a new tile. Therefore, no new columns appear in $T_{i+1}$ compared to $T_i$
and the columns of $T_{i+1}$ must still be linearly independent.

In the second case, the boxed entry belongs to a tile that is not already accounted for in the columns of $T_i$.
If this entry is contained in a fixed tile, we either have the same behavior as in the first case,
or it is a member of a free tile that starts in row $i+1$.
In the latter case, $T_{i+1}$ contains a new column, corresponding to this free tile.
The columns in $T_{i+1}$ are linearly independent in this case also.
\medskip

By the induction principle, we can conclude that all columns in $T$ are linearly independent,
so $G$ is a vertex of $\polyP_{\lambdavec/\muvec,\onevec}$.
\end{proof}

\bigskip

The following definition and observations allows us to reduce a number of cases to consider in a later argument.

Let $\lambdavec/\muvec$ and $\nuvec/\tauvec$ be two diagram shapes.
The \defin{disjoint union} of these diagrams is the shape
\[
 \lambdavec/\muvec \cup \nuvec/\tauvec = (\nu_1 + \lambdavec, \nuvec)/(\nu_1 + \muvec, \tauvec).
\]
For example, $(3,2)/(1) \cup (4,2,2)/(1,1)$ is the diagram
\[
(7,6,4,2,2)/(5,4,1,1) = 
 \young(:::::\hfil\hfil,::::\hfil\hfil,:\hfil\hfil\hfil,:\hfil,\hfil\hfil).
\]
Note that inserting $r$ empty rows and $c$ empty columns between the shapes
$\lambdavec/\muvec$ and $\nuvec/\tauvec$ in their disjoint union, a diagram of shape
\[
 (c+\nu_1 + \lambdavec, \underbrace{\nu_1,\dots,\nu_1}_r, \nuvec)/(c+\nu_1 + \muvec, \underbrace{\nu_1,\dots,\nu_1}_r, \tauvec)
\]
is obtained. GT-patterns of such a shape look like
\[
\begin{matrix}
\lambda'_1 &&  \dots &&  \lambda'_n && \nu_1 ,\dots,  \nu_1 && \nu_1 &&  \dots &&  \nu_l \\
& \ddots &&  && \ddots  && \ddots && \ddots &&   &&  \ddots  \\
&& \mu'_1 &&  \dots &&  \mu'_n && \nu_1 ,\dots, \nu_1 && \tau_1 &&  \dots && \tau_l \\
\end{matrix}
\]
where $\lambdavec' = c+\nu_1 + \lambdavec$ and  $\muvec' = c+\nu_1 + \muvec$.
Such GT-patterns consist of three distinct blocks: the first $n$ diagonals,
the next $r$ diagonals and the final $l$ diagonals.
Since $\mu'_n \geq \nu_1$, the south-west to north-east inequalities involving
the middle $r$ diagonals containing $\nu_1$ does not have any effect ---
for any choice of $r$ (even $r=0$), the valid fillings for different choices of $r$
are in natural correspondence with one another,
by inserting or deleting a number of diagonals from the middle block.

Likewise, there is a correspondence between patterns with different values of $c$,
by simply adding some constant to all elements in the first block.

Note that none of these two operations of changing $r$ or $c$ affect the weight of the pattern,
and it is fairly straightforward to see that changing $r$ or $c$ (\emph{i.e.} deleting or inserting
empty rows or columns between the shapes $\lambdavec/\muvec$ and $\nuvec/\tauvec$)
can be realized as invertible linear mappings between the corresponding GT-polytopes
that preserve lattice points.
In particular, this means that (integer) vertices are mapped to (integer) vertices.

Finally, note that the diagrams $(\lambdavec/\muvec) \cup (\nuvec/\tauvec)$ and 
$(\nuvec/\tauvec) \cup (\lambdavec/\muvec) $ are different in general.
However, there is an invertible, lattice-point-preserving linear map
between GT-patterns of these two shapes, namely the map that sends the pattern above to
\[
\begin{matrix}
\nu'_1 &&  \dots &&  \nu'_l && \lambda_1 ,\dots,  \lambda_1 && \lambda_1 &&  \dots &&  \lambda_n \\
& \ddots &&  && \ddots  && \ddots && \ddots &&   &&  \ddots  \\
&& \tau'_1 &&  \dots &&  \tau'_l && \lambda_1 ,\dots, \lambda_1 && \mu_1 &&  \dots && \mu_n \\
\end{matrix}
\]
with $\nuvec' = c+\lambda_1 + \nuvec$ and  $\tauvec' = c+\lambda_1 + \tauvec$.
This implies the following:
\begin{lemma}\label{lem:unionintegral}
The polytope $\polyP_{(\lambdavec/\muvec) \cup (\nuvec/\tauvec),\wvec}$ is integral if and
only if $\polyP_{(\nuvec/\tauvec)\cup (\lambdavec/\muvec),\wvec}$ is integral. 
\end{lemma}

\subsection{GT-polytopes corresponding to standard Young tableaux}

The integral GT-patterns in $\polyP_{\lambdavec/\muvec,\onevec}$
correspond to \defin{standard Young tableaux} with skew shape $\lambdavec/\muvec$,
that is, all boxes in the corresponding Young tableaux are different.

In this section, we completely characterize for which shapes $\lambdavec/\muvec$
the polytope $\polyP_{\lambdavec/\muvec,\onevec}$ is integral.
Because of the observation in the previous section, it is enough to consider shapes
$\lambdavec/\muvec$ with no empty row or column.

Given a partition $\lambdavec$, let $\lambdavec^+$ denote any partition
obtained from the diagram $\lambdavec$ by adding one box.
That is, $\lambdavec^+ \supset \lambdavec$ and $|\lambdavec^+| -|\lambdavec| = 1$.
Similarly, let $\lambdavec_-$ denote a partition obtained from the diagram
$\lambdavec$ by removing one box.

\begin{lemma}\label{lem:boxaddnonintegrality}
If $\polyP_{\lambdavec/\muvec,\onevec}$ is non-integral, then
$\polyP_{\lambdavec^+/\muvec,\onevec}$ and $\polyP_{\lambdavec/\muvec_-,\onevec}$
are also non-integral.
\end{lemma}
\begin{proof}
Let $G$ be a non-integer vertex in $\polyP_{\lambdavec/\muvec,\onevec}$:
\[
\begin{matrix}
\lambda_1 && \lambda_2 && \dots && \lambda_n \\
& \ddots && \ddots && && \ddots \\
&& \mu_1 && \mu_2 && \dots && \mu_n \\
\end{matrix}.
\]
The tiling matrix $T_G$ has full rank, since $G$ is a vertex.
Let $i$ be the row in the diagram $\lambdavec^+/\muvec$ where an extra box was added; consider
\[
G^+ = 
\begin{matrix}
\lambda_1 && \lambda_2 && \dots && \lambda_i+1 && \dots && \lambda_n \\
&\lambda_1 && \lambda_2 && \dots && \lambda_i && \dots && \lambda_n \\
&& \ddots && \ddots && && \ddots && && \ddots\\
&&& \mu_1 && \mu_2 && \dots && \mu_i && \dots  &&  \mu_n \\
\end{matrix}
\]
which is a point in $\polyP_{\lambdavec^+/\muvec,\onevec}$.
All $\lambda_j$ for $i \neq j$ in the second row from the top belong to a fixed tile of $G^+$.
Two things can happen:

\textbf{The entry $\lambda_i$ belongs to a fixed tile.}
This implies that the tiling matrix of $G^+$ is identical with that of $G$,
with the additional top row repeated twice.
This operation does not change the rank, so $G^+$ is a vertex.

\textbf{The entry $\lambda_i$ belongs to a free tile.} This free tile
was not present in $G$, so $T_{G^+}$ has an extra column.
However, this is the only free tile that intersects the second row from the top in $G^+$,
which implies that the corresponding column in the tiling matrix is linearly independent from the other columns.

It follows that $G^+$ is a non-integral vertex.
A similar argument shows that $G_-$ obtained from $G$ as
\[
G_- =
\begin{matrix}
\lambda_1 && \lambda_2 && \dots && \lambda_i && \dots && \lambda_n \\
& \ddots && \ddots && && \ddots && && \ddots\\
&& \mu_1 && \mu_2 && \dots && \mu_i && \dots  &&  \mu_n \\
&&& \mu_1 && \mu_2 && \dots && \mu_i-1 && \dots  &&  \mu_n \\
\end{matrix}
\]
is also a non-integral vertex.
\end{proof}

\begin{lemma}\label{lem:6shapes}
Whenever $\lambdavec/\muvec$ is any of the shapes
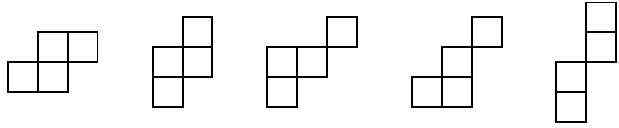
\[
\young(:\hfil\hfil,\hfil\hfil) \qquad  \young(:\hfil,\hfil\hfil,\hfil) \qquad  \young(::\hfil,\hfil\hfil,\hfil) \qquad  \young(::\hfil,:\hfil,\hfil\hfil)
\qquad  \young(:\hfil,:\hfil,\hfil,\hfil)
\]
the polytope $\polyP_{\lambdavec/\muvec,\onevec}$ is non-integral.
\end{lemma}
\begin{proof}
%
The shapes admit the following vertices:
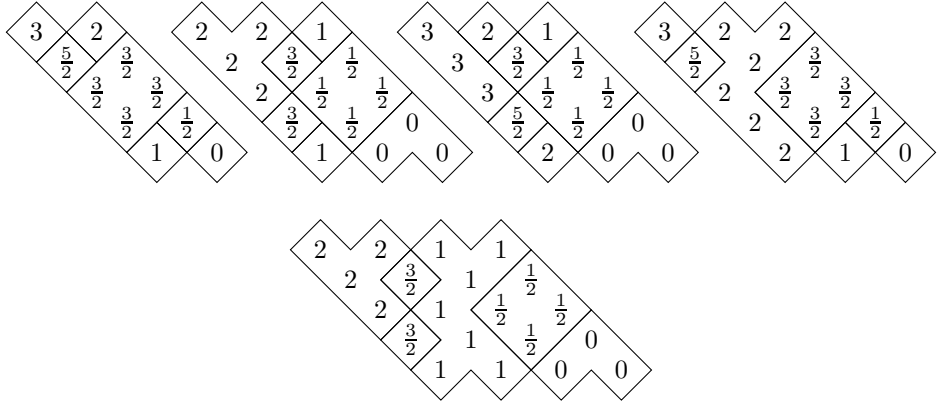
\begin{align*}
&
\begin{tikzpicture}[scale=0.4]
\draw[black] (2,0)--(3,-1)--(4,0)--(3,1)--(2,0)--cycle;
\draw[black] (4,0)--(5,-1)--(6,0)--(5,1)--(4,0)--cycle;
\draw[black] (6,-4)--(7,-5)--(8,-4)--(7,-3)--(6,-4)--cycle;
\draw[black] (8,-4)--(9,-5)--(10,-4)--(9,-3)--(8,-4)--cycle;
\draw[black] (3,-1)--(4,-2)--(5,-1)--(4,0)--(3,-1)--cycle;
\draw[black] (7,-3)--(8,-2)--(7,-1)--(6,0)--(5,-1)--(4,-2)--(5,-3)--(6,-4)--(7,-3)--cycle;
\draw[black] (7,-3)--(8,-4)--(9,-3)--(8,-2)--(7,-3)--cycle;
\node at (3,0) {$3$};
\node at (5,0) {$2$};
\node at (4,-1) {$\frac{5}{2}$};
\node at (6,-1) {$\frac{3}{2}$};
\node at (5,-2) {$\frac{3}{2}$};
\node at (7,-2) {$\frac{3}{2}$};
\node at (6,-3) {$\frac{3}{2}$};
\node at (8,-3) {$\frac{1}{2}$};
\node at (7,-4) {$1$};
\node at (9,-4) {$0$};
\end{tikzpicture} 
\mkern-52mu
\begin{tikzpicture}[scale=0.4]
\draw[black] (6,0)--(7,-1)--(8,0)--(7,1)--(6,0)--cycle;
\draw[black] (6,-4)--(7,-5)--(8,-4)--(7,-3)--(6,-4)--cycle;
\draw[black] (9,-3)--(8,-4)--(9,-5)--(10,-4)--(11,-5)--(12,-4)--(11,-3)--(10,-2)--(9,-3)--cycle;
\draw[black] (6,-2)--(5,-1)--(6,0)--(5,1)--(4,0)--(3,1)--(2,0)--(3,-1)--(4,-2)--(5,-3)--(6,-2)--cycle;
\draw[black] (5,-1)--(6,-2)--(7,-1)--(6,0)--(5,-1)--cycle;
\draw[black] (9,-3)--(10,-2)--(9,-1)--(8,0)--(7,-1)--(6,-2)--(7,-3)--(8,-4)--(9,-3)--cycle;
\draw[black] (5,-3)--(6,-4)--(7,-3)--(6,-2)--(5,-3)--cycle;
\node at (3,0) {$2$};
\node at (5,0) {$2$};
\node at (7,0) {$1$};
\node at (4,-1) {$2$};
\node at (6,-1) {$\frac{3}{2}$};
\node at (8,-1) {$\frac{1}{2}$};
\node at (5,-2) {$2$};
\node at (7,-2) {$\frac{1}{2}$};
\node at (9,-2) {$\frac{1}{2}$};
\node at (6,-3) {$\frac{3}{2}$};
\node at (8,-3) {$\frac{1}{2}$};
\node at (10,-3) {$0$};
\node at (7,-4) {$1$};
\node at (9,-4) {$0$};
\node at (11,-4) {$0$};
\end{tikzpicture}
\mkern-52mu
\begin{tikzpicture}[scale=0.4]
\draw[black] (4,0)--(5,-1)--(6,0)--(5,1)--(4,0)--cycle;
\draw[black] (6,0)--(7,-1)--(8,0)--(7,1)--(6,0)--cycle;
\draw[black] (6,-4)--(7,-5)--(8,-4)--(7,-3)--(6,-4)--cycle;
\draw[black] (6,-2)--(5,-1)--(4,0)--(3,1)--(2,0)--(3,-1)--(4,-2)--(5,-3)--(6,-2)--cycle;
\draw[black] (9,-3)--(8,-4)--(9,-5)--(10,-4)--(11,-5)--(12,-4)--(11,-3)--(10,-2)--(9,-3)--cycle;
\draw[black] (5,-1)--(6,-2)--(7,-1)--(6,0)--(5,-1)--cycle;
\draw[black] (9,-3)--(10,-2)--(9,-1)--(8,0)--(7,-1)--(6,-2)--(7,-3)--(8,-4)--(9,-3)--cycle;
\draw[black] (5,-3)--(6,-4)--(7,-3)--(6,-2)--(5,-3)--cycle;
\node at (3,0) {$3$};
\node at (5,0) {$2$};
\node at (7,0) {$1$};
\node at (4,-1) {$3$};
\node at (6,-1) {$\frac{3}{2}$};
\node at (8,-1) {$\frac{1}{2}$};
\node at (5,-2) {$3$};
\node at (7,-2) {$\frac{1}{2}$};
\node at (9,-2) {$\frac{1}{2}$};
\node at (6,-3) {$\frac{5}{2}$};
\node at (8,-3) {$\frac{1}{2}$};
\node at (10,-3) {$0$};
\node at (7,-4) {$2$};
\node at (9,-4) {$0$};
\node at (11,-4) {$0$};
\end{tikzpicture}
&
\mkern-62mu
\begin{tikzpicture}[scale=0.4]
\draw[black] (2,0)--(3,-1)--(4,0)--(3,1)--(2,0)--cycle;
\draw[black] (8,-4)--(9,-5)--(10,-4)--(9,-3)--(8,-4)--cycle;
\draw[black] (10,-4)--(11,-5)--(12,-4)--(11,-3)--(10,-4)--cycle;
\draw[black] (8,0)--(7,1)--(6,0)--(5,1)--(4,0)--(5,-1)--(4,-2)--(5,-3)--(6,-4)--(7,-5)--(8,-4)--(7,-3)--(6,-2)--(7,-1)--(8,0)--cycle;
\draw[black] (3,-1)--(4,-2)--(5,-1)--(4,0)--(3,-1)--cycle;
\draw[black] (9,-3)--(10,-2)--(9,-1)--(8,0)--(7,-1)--(6,-2)--(7,-3)--(8,-4)--(9,-3)--cycle;
\draw[black] (9,-3)--(10,-4)--(11,-3)--(10,-2)--(9,-3)--cycle;
\node at (3,0) {$3$};
\node at (5,0) {$2$};
\node at (7,0) {$2$};
\node at (4,-1) {$\frac{5}{2}$};
\node at (6,-1) {$2$};
\node at (8,-1) {$\frac{3}{2}$};
\node at (5,-2) {$2$};
\node at (7,-2) {$\frac{3}{2}$};
\node at (9,-2) {$\frac{3}{2}$};
\node at (6,-3) {$2$};
\node at (8,-3) {$\frac{3}{2}$};
\node at (10,-3) {$\frac{1}{2}$};
\node at (7,-4) {$2$};
\node at (9,-4) {$1$};
\node at (11,-4) {$0$};
\end{tikzpicture}
\end{align*}
\[
 \begin{tikzpicture}[scale=0.4]
\draw[black] (11,-3)--(10,-4)--(11,-5)--(12,-4)--(13,-5)--(14,-4)--(13,-3)--(12,-2)--(11,-3)--cycle;
\draw[black] (6,-2)--(5,-1)--(6,0)--(5,1)--(4,0)--(3,1)--(2,0)--(3,-1)--(4,-2)--(5,-3)--(6,-2)--cycle;
\draw[black] (10,0)--(9,1)--(8,0)--(7,1)--(6,0)--(7,-1)--(6,-2)--(7,-3)--(6,-4)--(7,-5)--(8,-4)--(9,-5)--(10,-4)--(9,-3)--(8,-2)--(9,-1)--(10,0)--cycle;
\draw[black] (5,-1)--(6,-2)--(7,-1)--(6,0)--(5,-1)--cycle;
\draw[black] (11,-3)--(12,-2)--(11,-1)--(10,0)--(9,-1)--(8,-2)--(9,-3)--(10,-4)--(11,-3)--cycle;
\draw[black] (5,-3)--(6,-4)--(7,-3)--(6,-2)--(5,-3)--cycle;
\node at (3,0) {$2$};
\node at (5,0) {$2$};
\node at (7,0) {$1$};
\node at (9,0) {$1$};
\node at (4,-1) {$2$};
\node at (6,-1) {$\frac{3}{2}$};
\node at (8,-1) {$1$};
\node at (10,-1) {$\frac{1}{2}$};
\node at (5,-2) {$2$};
\node at (7,-2) {$1$};
\node at (9,-2) {$\frac{1}{2}$};
\node at (11,-2) {$\frac{1}{2}$};
\node at (6,-3) {$\frac{3}{2}$};
\node at (8,-3) {$1$};
\node at (10,-3) {$\frac{1}{2}$};
\node at (12,-3) {$0$};
\node at (7,-4) {$1$};
\node at (9,-4) {$1$};
\node at (11,-4) {$0$};
\node at (13,-4) {$0$};
\end{tikzpicture}
\]
\end{proof}

\begin{lemma}\label{lem:3colshape}
Whenever $\lambdavec/\muvec$ is of the shape
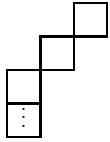
\[
\ytableausetup{centertableaux,boxsize=1.2em}
\begin{ytableau}
\none & \none &  \\
\none &  \\
\; \\
\scriptstyle{\vdots} \\
\end{ytableau}
\]
and the total number of boxes is at least four, the polytope $\polyP_{\lambdavec/\muvec,\onevec}$ is non-integral.
\end{lemma}
\begin{proof}
The cases $k=4$ and $k=5$ permit the following non-integral vertices:
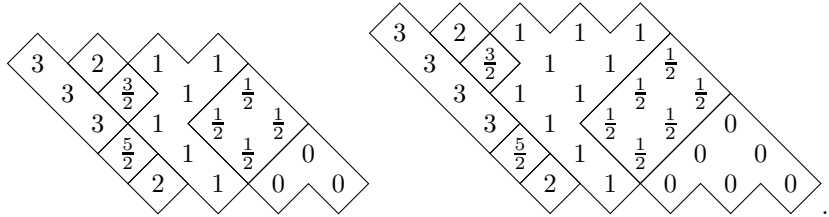
\[
 \begin{tikzpicture}[scale=0.4]
\draw[black] (4,0)--(5,-1)--(6,0)--(5,1)--(4,0)--cycle;
\draw[black] (6,-4)--(7,-5)--(8,-4)--(7,-3)--(6,-4)--cycle;
\draw[black] (6,-2)--(5,-1)--(4,0)--(3,1)--(2,0)--(3,-1)--(4,-2)--(5,-3)--(6,-2)--cycle;
\draw[black] (11,-3)--(10,-4)--(11,-5)--(12,-4)--(13,-5)--(14,-4)--(13,-3)--(12,-2)--(11,-3)--cycle;
\draw[black] (10,0)--(9,1)--(8,0)--(7,1)--(6,0)--(7,-1)--(6,-2)--(7,-3)--(8,-4)--(9,-5)--(10,-4)--(9,-3)--(8,-2)--(9,-1)--(10,0)--cycle;
\draw[black] (5,-1)--(6,-2)--(7,-1)--(6,0)--(5,-1)--cycle;
\draw[black] (11,-3)--(12,-2)--(11,-1)--(10,0)--(9,-1)--(8,-2)--(9,-3)--(10,-4)--(11,-3)--cycle;
\draw[black] (5,-3)--(6,-4)--(7,-3)--(6,-2)--(5,-3)--cycle;
\node at (3,0) {$3$};
\node at (5,0) {$2$};
\node at (7,0) {$1$};
\node at (9,0) {$1$};
\node at (4,-1) {$3$};
\node at (6,-1) {$\frac{3}{2}$};
\node at (8,-1) {$1$};
\node at (10,-1) {$\frac{1}{2}$};
\node at (5,-2) {$3$};
\node at (7,-2) {$1$};
\node at (9,-2) {$\frac{1}{2}$};
\node at (11,-2) {$\frac{1}{2}$};
\node at (6,-3) {$\frac{5}{2}$};
\node at (8,-3) {$1$};
\node at (10,-3) {$\frac{1}{2}$};
\node at (12,-3) {$0$};
\node at (7,-4) {$2$};
\node at (9,-4) {$1$};
\node at (11,-4) {$0$};
\node at (13,-4) {$0$};
\end{tikzpicture}
\begin{tikzpicture}[scale=0.4]
\draw[black] (4,0)--(5,-1)--(6,0)--(5,1)--(4,0)--cycle;
\draw[black] (7,-5)--(8,-6)--(9,-5)--(8,-4)--(7,-5)--cycle;
\draw[black] (6,-4)--(7,-3)--(6,-2)--(5,-1)--(4,0)--(3,1)--(2,0)--(3,-1)--(4,-2)--(5,-3)--(6,-4)--cycle;
\draw[black] (13,-5)--(14,-6)--(15,-5)--(16,-6)--(17,-5)--(16,-4)--(15,-3)--(14,-2)--(13,-3)--(12,-4)--(11,-5)--(12,-6)--(13,-5)--cycle;
\draw[black] (12,0)--(11,1)--(10,0)--(9,1)--(8,0)--(7,1)--(6,0)--(7,-1)--(6,-2)--(7,-3)--(8,-4)--(9,-5)--(10,-6)--(11,-5)--(10,-4)--(9,-3)--(10,-2)--(11,-1)--(12,0)--cycle;
\draw[black] (5,-1)--(6,-2)--(7,-1)--(6,0)--(5,-1)--cycle;
\draw[black] (12,-4)--(13,-3)--(14,-2)--(13,-1)--(12,0)--(11,-1)--(10,-2)--(9,-3)--(10,-4)--(11,-5)--(12,-4)--cycle;
\draw[black] (6,-4)--(7,-5)--(8,-4)--(7,-3)--(6,-4)--cycle;
\node at (3,0) {$3$};
\node at (5,0) {$2$};
\node at (7,0) {$1$};
\node at (9,0) {$1$};
\node at (11,0) {$1$};
\node at (4,-1) {$3$};
\node at (6,-1) {$\frac{3}{2}$};
\node at (8,-1) {$1$};
\node at (10,-1) {$1$};
\node at (12,-1) {$\frac{1}{2}$};
\node at (5,-2) {$3$};
\node at (7,-2) {$1$};
\node at (9,-2) {$1$};
\node at (11,-2) {$\frac{1}{2}$};
\node at (13,-2) {$\frac{1}{2}$};
\node at (6,-3) {$3$};
\node at (8,-3) {$1$};
\node at (10,-3) {$\frac{1}{2}$};
\node at (12,-3) {$\frac{1}{2}$};
\node at (14,-3) {$0$};
\node at (7,-4) {$\frac{5}{2}$};
\node at (9,-4) {$1$};
\node at (11,-4) {$\frac{1}{2}$};
\node at (13,-4) {$0$};
\node at (15,-4) {$0$};
\node at (8,-5) {$2$};
\node at (10,-5) {$1$};
\node at (12,-5) {$0$};
\node at (14,-5) {$0$};
\node at (16,-5) {$0$};
\end{tikzpicture}.
\]
These patterns can be generalized to any $k \geq 4$, by considering the patterns
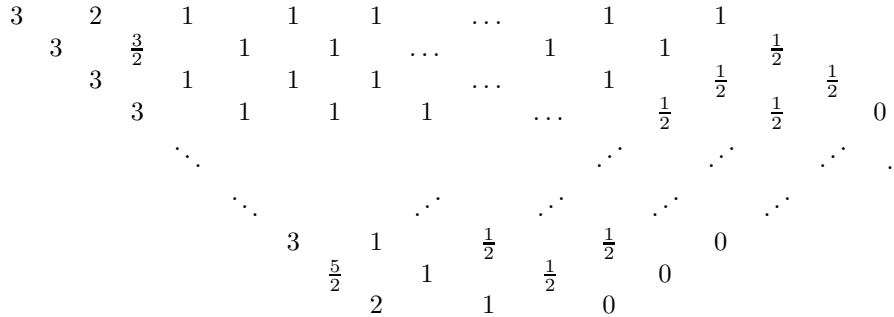
\[
\setcounter{MaxMatrixCols}{25}
\begin{matrix}
3 && 2 && 1 && 1 && 1 && \dots && 1 && 1\\
& 3 && \frac32 && 1 && 1 && \dots && 1 && 1 && \frac12 \\
&& 3 && 1 && 1 && 1 && \dots && 1 && \frac12 && \frac12 \\
&&& 3 && 1 && 1 && 1 && \dots && \frac12 && \frac12 && 0\\
&&&& \ddots &&  &&  &&  && \iddots && \iddots && \iddots  \\
&&&&& \ddots &&  && \iddots && \iddots && \iddots && \iddots  \\
&&&&&& 3 && 1 &&\frac12 && \frac12 && 0  \\
&&&&&&& \frac52 && 1 &&\frac12 && 0  \\
&&&&&&&& 2 && 1 && 0 
\end{matrix}.
\]
That the rows satisfy all the inequalities and that the row-sums increase by one upwards
is straightforward to check. Also, all integers in the pattern belong to fixed tiles
and the tiles containing $\tfrac12$, $\tfrac32$ and $\tfrac52$ produce a tiling matrix with full rank.
Thus, this is a non-integer vertex for every $k \geq 2$.
\end{proof}

Using \cref{lem:unionintegral} and \cref{lem:boxaddnonintegrality} together with \cref{lem:6shapes} and \cref{lem:3colshape}
it is clear that a large number of polytopes $\polyP_{\lambdavec/\muvec,\onevec}$ are non-integral.
In particular, all shapes in \eqref{eq:forbidden1} and \eqref{eq:forbidden2} below
have non-integral $\polyP_{\lambdavec/\muvec,\onevec}$ and adding boxes to these shapes preserves non-integrality.
The next lemma characterizes the diagrams that \emph{cannot} be obtained in this fashion.

Some terminology: A diagram $D_1$ is a \defin{subdiagram} of $D_2$,
if $D_1$ can be obtained from $D_2$ by removal of some boxes and deletion of empty rows and columns.

\begin{lemma} \label{lem:integralshapes}
A diagram without empty rows or columns that does not contain any of the diagrams
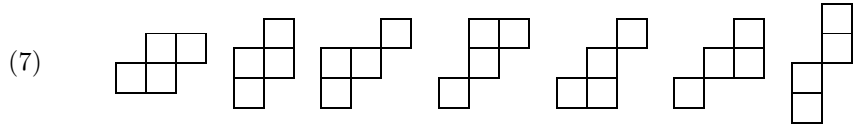
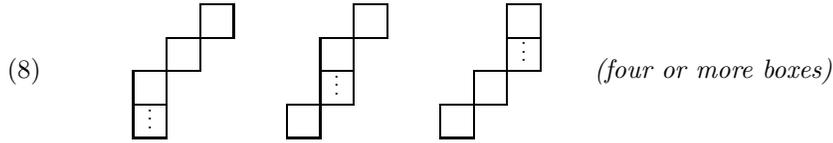
\begin{equation} \label{eq:forbidden1}
\young(:\hfil\hfil,\hfil\hfil) \quad  \young(:\hfil,\hfil\hfil,\hfil) \quad  \young(::\hfil,\hfil\hfil,\hfil)
\quad \young(:\hfil\hfil,:\hfil,\hfil) \quad  \young(::\hfil,:\hfil,\hfil\hfil) \quad  \young(::\hfil,:\hfil\hfil,\hfil)
\quad  \young(:\hfil,:\hfil,\hfil,\hfil)
\end{equation}
\begin{equation} \label{eq:forbidden2}
\ytableausetup{centertableaux,boxsize=1.2em}
\begin{ytableau}
\none & \none &  \\
\none &  \\
\; \\
\scriptstyle{\vdots} \\
\end{ytableau}
\qquad
\begin{ytableau}
\none & \none &  \\
\none &  \\
\none & \scriptstyle{\vdots} \\
\; \\
\end{ytableau}
\qquad
\begin{ytableau}
\none & \none &  \\
\none &  \none & \scriptstyle{\vdots}\\
\none & \; \\
\; \\
\end{ytableau} \qquad \text{(four or more boxes)}
\end{equation}
as a subdiagram is either a \emph{disjoint union of rows of boxes}, or of one of the shapes
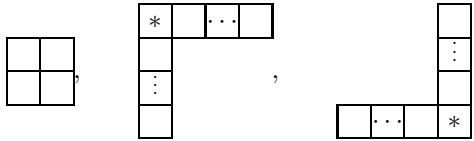
\[
\ytableausetup{centertableaux,boxsize=1.2em}
\begin{ytableau}
\; & \; \\ 
\; & \; \\
\end{ytableau}
,\qquad
\begin{ytableau}
\ast & & \cdots  & \\
   & \none\\
\scriptstyle{\vdots}  & \none\\
  & \none\\
\end{ytableau}
,\qquad
\begin{ytableau}
\none & \none & \none  &  \\
\none & \none & \none  & \scriptstyle{\vdots} \\
\none & \none & \none  &  \\
      & \cdots &  & \ast  \\
\end{ytableau}
\]
(called the $2\times2$-box, hook and reverse hook) where the box marked $\ast$ can either be present or not.
\end{lemma}
\begin{proof}
Assume we are given a diagram that avoids all the ten \emph{forbidden} diagrams in \eqref{eq:forbidden1} and \eqref{eq:forbidden2}.
A diagram that contains at most one box per column is a union of rows of boxes and it is easy to see that such
an arrangement does not contain a forbidden pattern.
Hence, assume that the diagram contains at least one instance where one box is on top of another
and consider the topmost such instance. If there are several such arrangements, pick the rightmost one.

Furthermore, we can assume that the diagram contains at least four boxes, since no arrangement of three boxes can be forbidden,
and it is easy to see that all such arrangements of at most three boxes is a hook, reverse hook or a union of rows.
Remember that we assume that there are no empty rows or columns, so there must be a third box
placed somewhere in one of the following ways: 
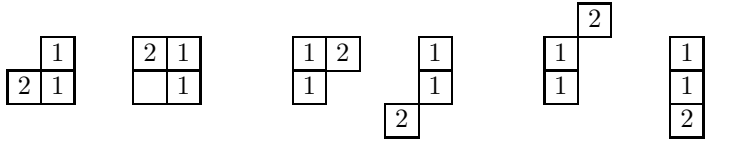
\[
\begin{ytableau}
\none & \none & \none \\
\none & 1    & \none \\
    2 & 1    & \none \\
\none & \none & \none \\
\end{ytableau}
\quad
\begin{ytableau}
\none & \none & \none \\
    2 & 1    & \none \\
      & 1    & \none \\
\none & \none & \none \\
\end{ytableau}
\quad
\begin{ytableau}
\none & \none & \none \\
\none & 1    &    2 \\
\none & 1    & \none \\
\none & \none & \none \\
\end{ytableau}
\quad
\begin{ytableau}
\none & \none & \none \\
\none & 1    & \none \\
\none & 1    & \none \\
   2  & \none & \none \\
\end{ytableau}
\quad
\begin{ytableau}
\none & \none & 2 \\
\none & 1    & \none \\
\none & 1    & \none \\
\none & \none & \none \\
\end{ytableau}
\quad
\begin{ytableau}
\none & \none & \none \\
\none & 1    & \none \\
\none & 1    & \none \\
\none & 2    & \none \\
\end{ytableau}.
\]
In the second case above, an extra box is indicated, since it must be present
in order for the shape to be a proper skew shape. In the last case, we can assume that the two boxes marked with $1$ are in the two top rows,
otherwise, one can consider the fifth case instead.

In the first five arrangements, some positions cannot have a box since that would introduce a forbidden pattern,
or would contradict the choice of the first two boxes. These are marked with $\times$ and we have the following possible arrangements:
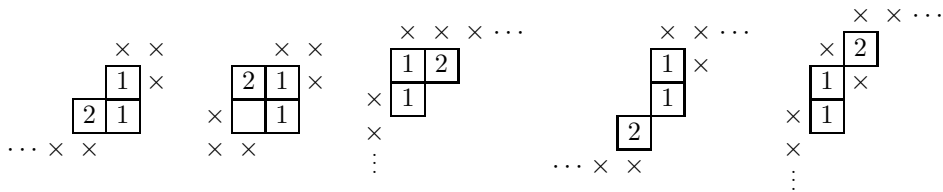
\[
\begin{ytableau}
\none &\none &\none & \none[\times]& \none[\times] \\
\none &\none &\none & 1     & \none[\times] \\
\none &\none &    2 & 1     & \none \\
\none[\cdots] &\none[\times] &\none[\times] & \none & \none \\
\end{ytableau}
\quad
\begin{ytableau}
\none &\none & \none[\times] & \none[\times] \\
\none &    2 & 1      & \none[\times] \\
\none[\times] &      & 1      & \none \\
\none[\times] & \none[\times]  & \none & \none\\
\end{ytableau}
\quad
\begin{ytableau}
\none  & \none[\times] & \none[\times] & \none[\times] & \none[\cdots] \\
\none  & 1    &    2 \\
\none[\times] & 1    & \none \\
\none[\times] & \none & \none \\
\none[\scriptstyle{\vdots}] & \none & \none \\
\end{ytableau}
\quad
\begin{ytableau}
\none &\none &\none &\none[\times] & \none[\times] & \none[\cdots]\\
\none &\none &\none & 1    & \none[\times] \\
\none &\none &\none & 1    & \none \\
\none &\none &   2  & \none & \none \\
\none[\cdots]& \none[\times]  & \none[\times]  & \none & \none \\
\end{ytableau}
\quad
\begin{ytableau}
\none & \none    & \none[\times] & \none[\times] & \none[\cdots] \\
\none & \none[\times] & 2 \\
\none & 1    & \none[\times] \\
\none[\times] & 1    & \none \\
\none[\times] & \none & \none \\
\none[\scriptstyle{\vdots}] & \none & \none \\
\end{ytableau}
\]
Knowing that some boxes are not present, other positions can be excluded as well, since the diagram must be a proper skew shape.
From here, it is straightforward to deduce that only the $2\times 2$-box or a (reverse) hook can be obtained by adding boxes in a non-forbidden fashion.

The last case is slightly different. If the diagram is only one column of boxes, it is a (degenerate) hook.
Since there are no empty rows or columns, we can assume the last case must be of one of the forms
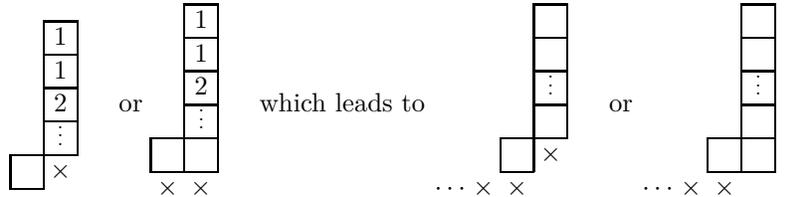
\[
\begin{ytableau}
\none & 1    & \none \\
\none & 1    & \none \\
\none & 2    & \none \\
\none & \scriptstyle{\vdots}    & \none \\
\; & \none[\times]  & \none \\
\end{ytableau}
\text{ or }
\begin{ytableau}
\none & 1    & \none \\
\none & 1    & \none \\
\none & 2    & \none \\
\none & \scriptstyle{\vdots}    & \none \\
\; &   & \none \\
\none[\times] & \none[\times]  & \none \\
\end{ytableau}
\text{ which leads to }
\begin{ytableau}
\none &\none &\none &        & \none \\
\none &\none &\none &        & \none \\
\none &\none &\none & \scriptstyle{\vdots}    & \none \\
\none &\none &\none &    & \none \\
\none &\none &    & \none[\times]  & \none \\
\none[\cdots] &\none[\times] & \none[\times]    & \none    & \none \\
\end{ytableau}
\text{ or }
\begin{ytableau}
\none &\none &\none &        & \none \\
\none &\none &\none &        & \none \\
\none &\none &\none & \scriptstyle{\vdots}    & \none \\
\none &\none &\none &    & \none \\
\none &\none &    &   & \none \\
\none[\cdots] &\none[\times] & \none[\times]    & \none    & \none \\
\end{ytableau}.
\]
It is now straightforward to deduce that these two cases must be reverse hooks.

Thus, every diagram that avoids the forbidden patterns is either a disjoin union of horizontal rows,
a $2\times 2$-box or some type of hook.
\end{proof}

It remains to show that the non-forbidden shapes in \cref{lem:integralshapes} are integral
to completely characterize all $\lambdavec/\muvec$ such that $\polyP_{\lambdavec/\muvec,\onevec}$ is integral.
We leave it as an exercise (this was also proved in \cite{King04stretched}) to show
that the $2 \times 2$-box gives a one-dimensional polytope, with two integer vertices
and proceed to show that whenever $\lambdavec/\muvec$ is a disjoint union of rows,
the corresponding Gelfand--Tsetlin polytope is integral.

\begin{proposition}\label{prop:integralribbon}
If $\lambdavec/\muvec$ is a disjoint union of rows, then $\polyP_{\lambdavec/\muvec,\wvec}$ is integral.
\end{proposition}
\begin{proof}
It is enough to consider the case when there are no empty rows or columns.
In that case, $\muvec = (\lambda_2,\lambda_3,\dots)$. Consider any GT-pattern with such a shape $\lambdavec/\muvec$:
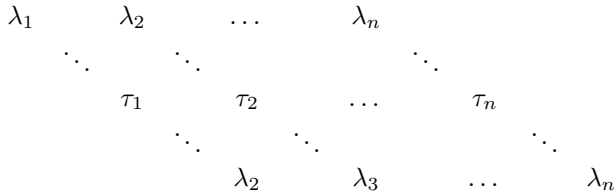
\[
\begin{matrix}
\lambda_1 && \lambda_2 && \dots && \lambda_n \\
& \ddots && \ddots && && \ddots \\
&& \tau_1 && \tau_2 && \dots && \tau_n \\
&&& \ddots && \ddots && && \ddots \\
&&&& \lambda_2 && \lambda_3 && \dots && \lambda_n
\end{matrix}
\]

Since $\lambda_j \geq \tau_j \geq \lambda_{j+1}$,
we can only have  $\tau_j = \tau_{j+1}$ if $\tau_j =  \lambda_{j+1} = \tau_{j+1}$.
Hence, $\tau_j$ and $\tau_{j+1}$ can only be members of the same tile,
if this tile also contains some element from the top and bottom row.

The conclusion is that no free tile can contain two (or more) entries from the same row.
Hence, all columns in the tiling matrix are of the form $(0,\dots,0,1\dots,1,0,\dots,0)$
and it is easy to show that such a matrix is totally unimodular.
Therefore, all points and especially all vertices in $\polyP_{\lambdavec/\muvec,\wvec}$
have totally unimodular tiling matrices.
Every minor of a totally unimodular matrix has determinant $0$ or $\pm 1$ and the result now follows from \cref{lem:unimodulartilingmatrix}.
\end{proof}

It now remains to show that hook shapes and reverse hook shapes give rise to integral polytopes as well.
This is done later in \cref{prop:integrallyclosedHooks}.

\section{Polytopes with the integer decomposition property and pulling triangulations}

The goal of this section is to show that every integral $\polyP_{\lambdavec/\muvec,\onevec}$
is a compressed polytope. We first recall some basic notions regarding
convex polytopes, see \cite{ziegler-lop-95,Haase14} for details.

Let $k\polyP$ denote the $k$-dilation of $\polyP$.
An integral polytope $\polyP \subset \setR^d$ is said to have the \defin{integer decomposition property} (IDP)
if for every $k \in \setN$ and $\xvec \in k \polyP \cap \setZ^d$,
we can find $\xvec^1,\xvec^2,\dots,\xvec^k \in \polyP \cap \setZ^d$ such that $\xvec^1 + \dots + \xvec^k = \xvec$.
The only simplices with the integer decomposition property are simplices with normalized volume one.
Such a simplex is called a \defin{unimodular simplex}.

In general, it is hard to determine if a polytope has the IDP using this definition,
but there are several stronger properties of polytopes that imply IDP.
For example, one can use the following proposition:
\begin{proposition}
If $\polyP$ has a triangulation into unimodular simplices, then $\polyP$ has the IDP.
Such a triangulation is called a \defin{unimodular triangulation}.
\end{proposition}

Information about the following definition can be found in \cite{Haase14,sullivant2006}:
\begin{definition}[Pulling triangulation]
Let $\polyP$ be a polytope and fix a total order on the vertices $p_1,\dots,p_k$ of $\polyP$.
The \defin{pulling triangulation} $\Delta_{pull}(\polyP)$ is defined recursively as follows:
If $p_1,\dots,p_k$ are affinely independent, then $\Delta_{pull}(\polyP)$ is just $\{\{p_1,\dots,p_k\}\}$.
Otherwise,
\[
\Delta_{pull}(\polyP) = \bigcup_F \{ \{p_k\} \cup \sigma | \sigma \in \Delta_{pull}(F) \},
\]
where the union is taken over all facets $F$ of $\polyP$, that do not contain $p_k$.
The ordering of the vertices in the facets are induced by the ordering on $\polyP$.
\end{definition}
One property of pulling triangulations, which follows easily from the definition, is the following:
If $\{p_0,\dots,p_d\} \in \Delta_{pull}(\polyP)$ (with the total order above),
then the minimal-dimensional face of $P$ that contains $\{p_0,\dots,p_j\}$ is $j$-dimensional,
for all $0\leq j \leq d$.
Finally, an integral convex polytope $\polyP$ is \defin{compressed} if every
pulling triangulation of $\polyP$ is a unimodular triangulation.

\bigskip

Consider an integral GT-pattern $G$ in $\polyP_{\lambdavec/\muvec,\onevec}$.
It is clear that $G$ is uniquely determined by the entries $(i,j)$ where $x^i_{j} > x^{i-1}_{j}$.
We call such an entry an \defin{increase}. 
Note that if $(i,j)$ is an increase, then $x^i_{j}$ and $x^{i-1}_{j}$ are members of different tiles.

\begin{lemma}\label{lem:tilingdimension}
Let $G_0, \dots, G_d$ be \emph{integer} vertices of $\polyP_{\lambdavec/\muvec,\onevec}$
such that the minimal-dimensional face of $\polyP_{\lambdavec/\muvec,\onevec}$ that contains $\{G_0, \dots, G_j\}$
is $j$-dimensional, for all $0 \leq j \leq d$. Then the simplex $\Delta = \{G_0, \dots, G_d \}$ is unimodular.
\end{lemma}
\begin{proof}
We use induction over $d$ and note that for $d=0$, the statement is clear.

Consider any integral $G$ in the $k$-dilation of $\Delta$, that is, any integral $G$ such that
\begin{equation}\label{eqn:coeffs}
G = a_0 G_0 + a_1 G_1 + \dots + a_d G_d
\end{equation}
where the $a_j \in \setR$ are non-negative numbers with sum $k$.
To prove that the simplex is unimodular, it suffices to show that the $a_j$ are integers.
It is enough to show that $a_d$ is must be an integer, since if $a_d$ is an integer,
\begin{equation}
G - a_d G_d = a_0 G_0 + a_1 G_1 + \dots + a_{d-1} G_{d-1}
\end{equation}
is an integer point in the $(k-a_d)$-dilation of the simplex $\Delta' = \{G_0,\dots, G_{d-1} \}$.
By induction, $\Delta'$ is unimodular, which implies that $a_j$ for $0\leq j \leq d-1$ are integers.
\medskip

We can now assume $a_d>0$ and observe that the tiling of $G$ must be a weak refinement of the tiling of $G_d$,
according to \cref{lem:tilingrefinement}.
Furthermore, we know that the minimal-dimensional face of $\polyP_{\lambdavec/\muvec,\onevec}$ that contains 
$\{G_0, \dots, G_{d-1}\}$ is $(d-1)$-dimensional and that the minimal-dimensional face
that contains $\{G_0, \dots, G_{d}\}$ is $d$-dimensional.
Hence, the tiling of $G_0+\dots+G_d$ is a \emph{strict} refinement of the tiling of $G_0+\dots+G_{d-1}$.
The tiling of $G_d$ is uniquely defined by the set of increases.
Thus, there must be an increase $(i,j)$ in $G_d$ which is not present in any of the $G_j$ for $j<d$,
and $x^i_{j} - x^{i-1}_{j}=1$ in $G_d$, since $G_d$ is an integer vertex.

From this, it follows that the difference $x^i_{j} - x^{i-1}_{j}$ in $G$ is given by $a_d$.
Since $G$ is an integral GT-pattern, $a_d$ is an integer.
\end{proof}

\begin{corollary}\label{onepartitioniscompressed}
If $\polyP_{\lambdavec/\muvec,\onevec}$ is an integer polytope, it is compressed.
\end{corollary}
\begin{proof}
Consider any simplex $\Delta = \{G_0,\dots,G_d\}$ in any pulling triangulation of $\polyP_{\lambdavec/\muvec,\onevec}$.
The vertices of $\Delta$ are also integer vertices of $\polyP_{\lambdavec/\muvec,\onevec}$, according to \cref{prop:onevertices}.
The properties of pulling triangulations ensure that the minimal-dimensional face of $\polyP_{\lambdavec/\muvec,\onevec}$ that contains $\{G_0, \dots, G_j\}$
is $j$-dimensional. \cref{lem:tilingdimension} now implies that $\Delta$ is unimodular.
Hence, every pulling triangulation of $\polyP_{\lambdavec/\muvec,\onevec}$ is unimodular, so the polytope is compressed.
\end{proof}

It is known that all lattice points in compressed polytopes are vertices of the polytope,
so this result cannot be extended to general $\wvec$.
For general $\wvec$, there are plenty of examples of integer GT-patterns that are not vertices.

\section{Conditional results on integrality and integer decomposition property}

In this section we establish integrality and IDP of $\polyP_{\lambdavec/\muvec,\wvec}$ under certain conditions.

\begin{lemma}\label{lem:connectedrowinsertion}
Let $G$ be the GT-pattern
\[
G=\begin{matrix}
\lambda_1 && \lambda_2 && \dots && \lambda_n \\
& \mu_1 && \mu_2 && \dots && \mu_n 
\end{matrix}
\]
with rational entries. Then for every $t \in \setQ$ between $|\muvec|$ and $|\lambdavec|$,
there is a $\nuvec$ such that $|\nuvec|=t$ and a GT-pattern
\[
G'=\begin{matrix}
\lambda_1 && \lambda_2 && \dots && \lambda_n \\
& \nu_1 && \nu_2 && \dots && \nu_n \\
&& \mu_1 && \mu_2 && \dots && \mu_n 
\end{matrix}
\]
such that the entries $\lambda_i$ and $\mu_j$ belong to the same
tile in $G'$ if $\lambda_i$ and $\mu_j$ belong to the same tile in $G$.
\end{lemma}
\begin{proof}
Let $\delta>0$ be the largest rational number such that $\frac{1}{\delta}G$
is an integer GT-pattern and $t/\delta$ is an integer. It suffices to show the following statement:

\textbf{Statement:} For every $t = a\delta$ with $a \in \setN$ such that $|\muvec|\leq t \leq |\lambdavec|$,
there is a $\nuvec$ such that $|\nuvec|=t$, $\frac{1}{\delta}\nuvec$ is an integer partition,
\[
G'=\begin{matrix}
\lambda_1 && \lambda_2 && \dots && \lambda_n \\
& \nu_1 && \nu_2 && \dots && \nu_n \\
&& \mu_1 && \mu_2 && \dots && \mu_n
\end{matrix}
\]
is a GT-pattern and the entries $\lambda_i$ and $\mu_j$ belong to the same
tile in $G'$ if $\lambda_i$ and $\mu_j$ belong to the same tile in $G$.
\medskip 

It is clear that the statement is true for $t=|\muvec|$ and $t=|\lambdavec|$
since then we may choose $\nuvec=\muvec$ and $\nuvec=\lambdavec$ respectively.
Using induction, it is enough to show that if the statement is
true for $t = a \delta$, then it is also true for $t = (a+1)\delta$.
\medskip 

Now consider $G'$ for some $\nuvec$ satisfying the conditions in the statement.
It suffices to find some $\nu_i$ that can be increased by $\delta$,
such that if $\mu_i=\lambda_{i+1}$ then at least one of the equalities $\mu_i=\nu_i$ or
$\mu_i=\nu_{i+1}$ holds.
Keep in mind that $\lambdavec \interl \muvec$ and that $\lambdavec \interl \nuvec \interl \muvec$.
Furthermore, $\delta$ is less than or equal to the smallest non-zero difference between entries in $G'$.
There are four cases to consider:
\begin{enumerate}
\item[Case 1:] There is some $i$ such that $\lambda_i>\nu_i>\mu_i$. 
In this case, $\nu_i$ does not belong to a tile in $G'$, so it can be increased by $\delta$ without breaking the
property of connected tiles.

\item[Case 2:] There is some $i<n$ such that $\lambda_i>\nu_i = \mu_i$ but $\nu_i>\lambda_{i+1}$.
Thus, $\nu_i$ and $\mu_i$ belong to the same tile, but $\nu_i$ does not belong to a tile intersecting the top row.
Hence, $\nu_i$ can be increased by $\delta$ without breaking the property of connected tiles.

\item[Case 3a:] For $1\leq j < i$ we have $\lambda_i=\nu_i=\mu_i$. Then the pattern $G'$ is of the form
\[
\begin{matrix}
\lambda_i && \nu_i && \nu_{i+1} && \nu_{i+2} && \dots && \nu_{i+l} && \dots\\
& \nu_i && \nu_{i+1}&& \nu_{i+2} && \dots && \underline{\nu_{i+l}} && \nu_{i+l} && \dots\\
&& \nu_i && \nu_{i+1}&& \nu_{i+2}  && \dots && \nu_{i+l} && \dots
\end{matrix}
\]
for some $l\geq 0$ and we can safely increase $\nu_{i+l}$ since increasing 
this will not break any tile into two smaller tiles.

\item[Case 3b:] The last case can essentially be considered as an extension of Case 2 or 3a,
\[
G'=
\begin{matrix}
\lambda_i && \nu_i && \nu_{i+1}  && \dots && \nu_{n-1}\\
& \nu_i && \nu_{i+1}&& \dots && \nu_{n-1} && \underline{\nu_{n}}\\
&& \nu_i && \nu_{i+1}&& \dots && \nu_{n-1} && \nu_{n}
\end{matrix}
\]
and it is clear that the underlined entry can be increased, since $\nu_{n-1}>\nu_n$.
\end{enumerate}
The cases above cover all possibilities and from these observations, the statement follows.
\end{proof}

\begin{example}\label{ex:rowinsertion}
We illustrate the previous lemma with an example. Let $G$ be given by
\[
G=\begin{matrix}
4 && \mathbf{\frac52} && \frac32 && \underline{0} \\
& \mathbf{\frac52} && 0 && \underline{0} && \underline{0}
\end{matrix}.
\]
Suppose we wish to insert a line in the middle, where the sum of the entries is $7$.
There are two tile that needs to be preserved, with content $\frac52$ and $\underline{0}$ respectively.
Start with the pattern
\[
\begin{matrix}
4 && \frac52 && \frac32 && 0 \\
& \frac52 && 0 && 0 && 0 \\
&& \frac52 && 0 && 0 && 0
\end{matrix}.
\]
Case 1 does not apply, but Case 2 in the lemma works. This gives
\[
\begin{matrix}
4 && \frac52 && \frac32 && 0 \\
& \frac52 && \frac52 && 0 && 0 \\
&& \frac52 && 0 && 0 && 0
\end{matrix}.
\]
The middle row sum is now 5, and we can now apply Case 2 again, followed by 3a:
\[
\begin{matrix}
4 && \frac52 && \frac32 && 0 \\
& \frac52 && \frac52 && \frac32 && 0 \\
&& \frac52 && 0 && 0 && 0
\end{matrix}
\quad
\longrightarrow 
\quad
\begin{matrix}
4 && \mathbf{\frac52} && \frac32 &&  \underline{0} \\
& 3 && \mathbf{\frac52} && \frac32 &&  \underline{0} \\
&& \mathbf{\frac52}&& 0 &&  \underline{0} &&  \underline{0}
\end{matrix}.
\]
Note that the leftmost 0 in the bottom row was \emph{not} part of a larger tile in $G$,
so it does not have to be part of the larger $\underline{0}$-tile. 
\end{example}

\begin{proposition}[Integrality and refinement]\label{prop:integralityrefinement}
Let  $\wvec' \refin \wvec$.
If $\polyP_{\lambdavec/\muvec,\wvec}$ is a non-integral polytope,
then $\polyP_{\lambdavec/\muvec,\wvec'}$ is also non-integral.
\end{proposition}
\begin{proof}
Let $\wvec$ be given by $(w_1,\dots,w_{i-1},w_i+w_i',w_{i+1},\dots,w_n)$
and let $\wvec'$ be $(w_1,\dots,w_{i-1},w_i,w_i',w_{i+1},\dots,w_n)$.
It is clear that $\wvec' \refin \wvec$ and that any other refinement
can be obtained by repeating this type of refinement.

Using Lemma \ref{lem:connectedrowinsertion}, any rational point $G$ in the polytope $\polyP_{\lambdavec/\muvec,\wvec}$ 
can be mapped to a point in $\polyP_{\lambdavec/\muvec,\wvec'}$ by inserting some new row
between rows $i-1$ and $i$ in the GT-pattern $G$.
Furthermore, every tile in $G$ can be naturally identified with a tile in $G'$.
Tiles in $G$ that do not cross from row $i-1$ to $i$ are preserved identically in $G'$
and each tile in $G$ that does cross this line also appear in $G'$ (possibly with some extra elements),
since Lemma \ref{lem:connectedrowinsertion} guarantees that no
tile in $G$ is broken into two or more tiles in $G'$.

\emph{The only new tiles that appear in $G'$ are tiles consisting of only one element in row $i$.}
Therefore, the tiling matrices of $G$ and $G'$ can be written as
\[
T_G =
\begin{pmatrix}
t_{11} & t_{12} & \dots & t_{1r} \\
\vdots & \vdots & \ddots & \vdots \\
t_{j1} & t_{i2} & \dots & t_{jr} \\
\vdots & \vdots & \ddots & \vdots \\
t_{m1} & t_{m2} & \dots & t_{mr}
\end{pmatrix}
\text{         and         }
T_{G'}=
\begin{pmatrix}
t_{11} & t_{12} & \dots & t_{1r}  & 0 & \dots & 0\\
\vdots & \vdots & \ddots & \vdots & \vdots  & \ddots & \vdots \\
t'_{j1} & t'_{i2} & \dots & t'_{jr} & 1  & \dots & 1 \\
t_{j1} & t_{i2} & \dots & t_{jr} & 0& \dots & 0\\
\vdots & \vdots & \ddots & \vdots & \vdots& \ddots & \vdots\\
t_{m1} & t_{m2} & \dots & t_{mr} & 0 & \dots & 0
\end{pmatrix},
\]
where $j = m-i+1$ and $m$ is the number of rows in $G$.

Assume now that $G$ is a non-integral vertex of $\polyP_{\lambdavec/\muvec,\wvec}$.
If $G'$ is a vertex, we are done. 
Otherwise, \cref{lem:tilingkernel} implies that there is an some
$\yvec \in \ker T_{G'}$ such that adding $y_i$ to the entries in tile $i$ in $G'$ for all $i=1,2,\dots$,
we obtain a vertex of $\polyP_{\lambdavec/\muvec,\wvec'}$.
Note that if $\yvec$ is in the kernel of $T_{G'}$, we must have $y_1=\dotsb=y_r=0$,
otherwise, $(y_1,\dots,y_r)$ would be a non-zero entry in
the kernel of $T_G$, which is impossible since $G$ is a vertex.
Hence, the vertex constructed by adding $y_i$ to the tiles of $G'$ preserves all entries in the tiles corresponding to the first $r$
columns in $T_{G'}$. Since there were non-integral entries among these, the vertex constructed in this manner is non-integral.
\end{proof}

The previous proposition is a bit technical, so an example is justified:
\begin{example}
Consider $\polyP_{\lambdavec,\wvec}$ for $\lambda=(4,4,2,1,0)$ and $\wvec = (1,2,2,3,3)$.
This polytope is non-integral, since it has the following GT-pattern as a vertex:
\[
G=\begin{matrix}
4 &&  4 && 2  && 1  && 0  \\
& 4 && \frac52 && \frac32 && 0  \\
 && \frac52 && \frac52 && 0\\
 &&& \frac52 && \frac12 \\
 &&&& 1
\end{matrix} \quad 
T_G = 
\begin{pmatrix}
0 & 0 & 0 \\
1 & 1 & 0 \\
2 & 0 & 0 \\
1 & 0 & 1 \\
0 & 0 & 0
\end{pmatrix}.
\]
We now consider the refinement where $\wvec' = (1,2,2,2,1,3)$.
Reusing the $\nuvec$ calculated in \cref{ex:rowinsertion},
the pattern $G'$ in $\polyP_{\lambdavec,\wvec'}$ is given as
\[
G'=\begin{matrix}
4 &&  4 && 2  && 1  && 0  && 0\\
& 4 && \frac52 && \frac32 && 0 && 0 \\
&& \underline{3} && \frac52 && \underline{\frac32} && 0  \\
 &&& \frac52 && \frac52 && 0\\
 &&&& \frac52 && \frac12 \\
 &&&&& 1
\end{matrix} \quad 
T_{G'} = 
\begin{pmatrix}
0 & 0 & 0 &0 & 0\\
1 & 1 & 0 &0 & 0 \\
1 & 0 & 0 & 1 & 1\\
2 & 0 & 0 &0 & 0\\
1 & 0 & 1 &0 & 0\\
0 & 0 & 0 &0 & 0
\end{pmatrix}.
\]
This GT-pattern is not a vertex, but we can use the vector $(0,0,0,1,-1)$
in the kernel of $T_{G'}$ to adjust the underlined tiles in $G'$:
\[
G''=\begin{matrix}
4 &&  4 && 2  && 1  && 0  && 0\\
& 4 && \frac52 && \frac32 && 0 && 0 \\
&& 4 && \frac52 && \frac12 && 0  \\
 &&& \frac52 && \frac52 && 0\\
 &&&& \frac52 && \frac12 \\
 &&&&& 1
\end{matrix}
 \quad 
T_{G''} = 
\begin{pmatrix}
0 & 0 & 0 &0 \\
1 & 1 & 0 &0  \\
1 & 0 & 0 & 1 \\
2 & 0 & 0 &0 \\
1 & 0 & 1 &0 \\
0 & 0 & 0 &0 
\end{pmatrix}.
\]
It is now clear that $G''$ is a vertex, so $\polyP_{\lambdavec,\wvec'}$ is non-integral.
\end{example}

\bigskip

The next proposition can be used as a main tool in an inductive argument to prove
IDP for general $\wvec$.
\begin{proposition}[Tableaux box refinement]\label{prop:boxcoloring}
Let  $\wvec' \refin \wvec$. 
If $\polyP_{\lambdavec/\muvec,\wvec'}$ has the IDP\footnote{Note that this assumption implies integrality of $\polyP_{\lambdavec/\muvec,\wvec'}$.},
then $\polyP_{\lambdavec/\muvec,\wvec}$ also has the IDP.
\end{proposition}
\begin{proof}
Let $\wvec = (w_1,\dots,w_{i-1},w_i+w_i',w_{i+1},\dots,w_n)$. It suffices to consider refinements
of the form $\wvec' = (w_1,\dotsc,w_i,w_i',w_{i+1},\dotsc,w_n)$.
Let $G$ be an integral GT-pattern in $\polyP^{k}_{\lambdavec/\muvec,\wvec}$. 
We need to show that $G$ can expressed as
\begin{equation}\label{eq:boxeq}
G = G_1 + G_2 + \dots + G_k,\quad G_i \in \Gset_{\lambdavec/\muvec,\wvec}.
\end{equation}
Consider the tableau $T$ that corresponds to $G$. 
There are $k(w_i+w'_i)$ boxes with content $i$ and no two of these boxes appears in the same column.
Hence, there is a natural ordering of these boxes, from left to right.
We now construct a new tableau $T'$ by first adding $1$ to all boxes with content greater than or equal to $i+1$,
followed by adding one to the $kw_1$ \emph{rightmost} boxes in $T$ with content $i$.
See \cref{fig:boxrefinement} for an example of this transformation.
\begin{figure}[!ht]
\[
T=
\young(:::111111223,:::222223334,223333334,334555), \quad T' = 
\young(:::111111224,:::222224445,223333445,335666)
\]
\caption{Example of the tableau transformation for the parameters $k=3$, $w=(2, 3, 4, 1, 1)$ and $w'=(2, 3, 2,2, 1, 1)$.}
\label{fig:boxrefinement}
\end{figure}
The tableau $T'$ represents an integral GT-pattern in $\polyP^{k}_{\lambdavec/\muvec,\wvec'}$.
Since $\polyP_{\lambdavec/\muvec,\wvec'}$ has the IDP,
there are tableaux $T_j'$ in $\Gset_{\lambdavec/\muvec,\wvec'}$
such that $T'=T'_1 \tinsert T'_2 \tinsert \dotsb \tinsert T'_k$.
By construction, any row in $T'$ that has a box with content $i$
cannot have a row with a box with content $i+1$ below it.
Since $T'$ is the $\tinsert$-sum of the $T'_j$'s, this property also holds for each individual $T'_j$.

Hence, the 
inverse transformation\footnote{All boxes with content $i+1$ are replaced with boxes with content $i$, 
and all boxes with content $\geq i$ have its content decreased by one.} 
can be applied to each $T'_j$ to obtain tableau $T_j$,
corresponding to elements in $\Gset_{\lambdavec/\muvec,\wvec}$.
Subsequently, $T = T_1\tinsert T_2 \tinsert \dotsb \tinsert T_k$, which is equivalent with \eqref{eq:boxeq}.
\end{proof}

\begin{proposition}\label{prop:integrallyclosedHooks}
All GT-polytopes $\polyP_{\lambdavec,\onevec}$
where $\lambdavec = (h,1,1,\dots,1)$, \emph{i.e.}, $\lambdavec$ is a hook, are integrally closed.
\end{proposition}
\begin{proof}
Let $T$ be tableau in $\polyP^k_{\lambdavec,\onevec}$.
It suffices to show that there is a $T'$ such that $T' \tinsert T''$
for some $T'$ in $ \polyP_{\lambdavec,\onevec}$ and $T''$ in $\polyP^{k-1}_{\lambdavec,\onevec}$.

Consider the first column in $T$, which is a subset of $\{1,\dots,|\lambdavec|\}$.
If some number $j$ is not present in this column,
then the $k$ boxes with content $j$ cannot be in the first $k$ columns of $T$
so there must be some other column consisting of a single box with content $j$.
For all such missing numbers $j$, we record a corresponding one-box column.
Let $T'$ be the first column, together with the recorded columns, in left-to-right order.
By construction, $T'$ is a standard tableau with shape $\lambdavec$ and the complement of these columns, $T''$,
is easily seen to be a semi-standard Young tableau in $\polyP^{k-1}_{\lambdavec,\onevec}$.

It is now straightforward to show that $T = T' \tinsert T''$.
\end{proof}
\begin{remark}
The same proof can now be carried out for reverse hooks, that is, skew shapes of the form
$\lambda = (h,h,\dots,h)$ with $\length(\lambda)=l$ and $\mu = (h-1,h-1,\dots,h-1)$ with $\length(\mu)=l-1$.
\end{remark}

\begin{corollary}\label{cor:completeIntegralitycharacterization}
Using \cref{prop:integralribbon} and \cref{prop:integrallyclosedHooks} and \cref{lem:integralshapes},
$\polyP_{\lambdavec/\muvec,\onevec}$ is integral if and only if $\lambdavec/\muvec$ is one of the shapes
\begin{itemize}
\item a disjoint union of rows,
\item a $2 \times 2$-box,
\item a hook, (possibly with the corner box missing).
\end{itemize}
\end{corollary}

\section{A partial order on GT-polytopes}

For fixed $\lambdavec/\muvec$, the partial order $\wvec' \refin \wvec$ 
induces an order on the polytopes $\polyP_{\lambdavec/\muvec,\wvec}$. 
This order has several nice properties given by
\cref{prop:integralityrefinement} and \cref{prop:boxcoloring}:

\begin{theorem}[Partial order properties]\label{thm:paoprop}
Let $\wvec' \refiner \wvec$
and let $P=\polyP_{\lambdavec/\muvec,\wvec} \subset \setR^d$ and $P'=\polyP_{\lambdavec/\muvec,\wvec'}\subset \setR^{d'}$.
Then
\begin{enumerate}
\item $|P' \cap \setZ^{d'}|$ is greater or equal to $|P \cap \setZ^{d}|$.
\item If $P'$ is empty, then $P$ is empty.
\item If $P'$ is integral, then $P$ is integral.
\item If $P'$ has the integer decomposition property, then so does $P$.
\end{enumerate}
\end{theorem}
The first item follows from using the same injection as described in \cref{prop:boxcoloring},
and the rest are evident from previous results.
We also conjecture one additional statement, supported by computer experiments:
\begin{enumerate}
\item[\textit{(5)}] If $P'$ is a unimodular simplex, then $P$ a unimodular simplex.
\end{enumerate}

In \cref{fig:polyorder}, \cref{thm:paoprop} is illustrated in the non-skew case $\lambdavec = (4,3,1)$.
The nodes are the values of $\wvec$ and arrows indicate the partial order $\refiner$.
Parameters $\wvec$ that give empty polytopes are not shown.
The solid and dashed frames indicate integral polytopes, whereas dotted frames indicate non-integral.
The solid frames correspond to unimodular simplices.

Note that it is possible to extend this figure where different permutations of $\wvec$ are presented.

\begin{figure}[!ht]
\centering
\begin{tikzpicture}[>=triangle 45,xscale=1.6,y=-1cm, yscale=1.2]
    \node[draw] (n422) at (1,0) {422};
    \node[draw] (n332) at (2,0) {332};
    \node[draw] (n431) at (3,0) {431};

    \node[draw,dashed] (n2222) at (0,1) {2222};
    \node[draw,dashed] (n3221) at (1,1) {3221};
    \node[draw] (n4211) at (2,1) {4211};
    \node[draw] (n3311) at (3,1) {3311};
    
\node[draw,dotted] (n22211) at (1,2) {22211};    
\node[draw,dashed] (n32111) at (2,2) {32111};
\node[draw] (n41111) at (3,2) {41111};

\node[draw,dotted] (n221111) at (2,3) {221111};
\node[draw,dotted] (n311111) at (3,3) {311111};

\node[draw,dotted] (n2111111) at (2,4) {2111111};

\node[draw,dotted] (n11111111) at (2,5) {11111111};

\draw [->] (n422) -- (n2222);
\draw [->] (n422) -- (n3221);
\draw [->] (n422) -- (n4211);
\draw [->] (n332) -- (n3221);
\draw [->] (n332) -- (n3311);
\draw [->] (n431) -- (n3221);
\draw [->] (n431) -- (n4211);
\draw [->] (n431) -- (n3311);
\draw [->] (n2222) -- (n22211);
\draw [->] (n3221) -- (n22211);
\draw [->] (n3221) -- (n32111);
\draw [->] (n4211) -- (n22211);
\draw [->] (n4211) -- (n32111);
\draw [->] (n4211) -- (n41111);
\draw [->] (n3311) -- (n32111);
\draw [->] (n22211) -- (n221111);
\draw [->] (n32111) -- (n221111);
\draw [->] (n32111) -- (n311111);
\draw [->] (n41111) -- (n311111);
\draw [->] (n41111) -- (n221111);
\draw [->] (n221111) -- (n2111111);
\draw [->] (n311111) -- (n2111111);
\draw [->] (n2111111) -- (n11111111);
\end{tikzpicture}
\caption{The partial order on the polytopes $\polyP_{431,\wvec}$ with different $\wvec$.}\label{fig:polyorder}
\end{figure}

\section{GT-polytopes without weight restriction}\label{sec:noweight}

For completeness, we briefly mention integrality and the integer decomposition property
in the case when no restriction on the weight is imposed.
This setting is significantly less complicated.

Given a skew shape $\lambdavec/\muvec$, define the convex polytope 
$\polyP_{\lambdavec/\muvec} \subset \setR^{mn}$
consisting of all GT-patterns $(x^i_j)_{1\leq i \leq m,1\leq j \leq n}$ that satisfy 
the equalities  $\xvec^1 = \lambdavec$ and $\xvec^m = \muvec$.
These polytopes can be seen as GT-polytopes without any restriction on the row sums.

\begin{proposition}
All $\polyP_{\lambdavec/\muvec}$ have the integer decomposition property. 
\end{proposition}
\begin{proof}[Proof sketch:]
A simple proof appears in \cite{Alexandersson20141} so an example illustrating the idea is enough.
Given $G \in \polyP^k_{\lambdavec/\muvec} \cap \setZ^{mn}$, consider the corresponding tableau.
\[
T=\young(::::::111115,:::111333,122222445,245) \text{ Here, } \lambdavec/\muvec=(4, 3, 3, 1)/(2, 1) \text{ and } k=3.
\]
Note that columns appear in blocks of $k$. By selecting the $j$th column in each block for $j=1,2,\dots,k$,
$k$ smaller tableaux are constructed and the big tableau can be expressed as the $\tinsert$-sum of the smaller 
tableaux (recall the definition of this operation in \cref{ssec:adding}).
In this particular case,
\[
T= 
\young(::11,:13,124,2)
\tinsert
\young(::11,:13,224,4)
\tinsert
\young(::15,:13,225,5)
\]
This construct shows that $\polyP_{\lambdavec/\muvec}$ has the IDP.
\end{proof}
Note that this implies that all $\polyP_{\lambdavec/\muvec}$  are also integral.
A recent paper \cite{GusevKT13} studies the number of vertices of this type of polytope.
There seem to be a lot of open questions in this area.

\section{Connection with contingency matrices}\label{sec:contingencyMatrices}

There is a natural correspondence between \emph{contingency matrices} and certain GT-polytopes.
A contingency matrix is a matrix of non-negative integers, with specified row sums and column sums.
Such matrices appear naturally in statistics as well as
in representation theory (see \emph{e.g.}~\cite{Diaconis1995}) and many other areas.

Let $\lambdavec/\muvec$ be a disjoint union of rows and
consider $G \in \polyP^k_{\lambdavec/\muvec,\onevec}$ and the
corresponding tableau.
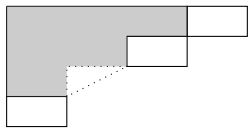
\begin{center}
\begin{tikzpicture}[scale=0.4,y=-1cm]
\foreach \point [count=\i] in
{(6,0),(4,1),(0,3)}
{
\begin{scope}[shift={\point}]
    \draw (0,0) --(2,0) --  (2,1) --(0,1) -- cycle;
\end{scope}
}
\begin{scope}[style=dotted]
\draw (4,2)--(2,2)--(2,3)--cycle;
\end{scope}
 \draw (0,0) -- (6,0) -- (6,1) -- (4,1);
 \draw (2,3)-- (0,3)--(0,0);
\fill[fill=gray,fill opacity=0.4]
(0,0) -- (6,0) -- (6,1) -- (4,1) -- (4,2)-- (2,2)-- (2,3)-- (0,3) -- cycle; 
\end{tikzpicture}
\end{center}

Note that we might need to generalize the notion of a Young tableau slightly,
where boxes can have fractional width.
Each row $i$ contains $(\lambda_i-\mu_i)k$ boxes and there are $m$ rows.
Let $a_{ij}$ be the number of (possibly fractional) boxes with content $j$ in row $i$.
The quantities $a_{ij}$ can be computed using the observation in \cref{rem:gtboxcount}.
Now consider the matrix
\[
A=
\begin{pmatrix}
a_{11} & a_{12} & \dots & a_{1n} \\
a_{21} & a_{22} & \dots & a_{2n} \\
\vdots & \vdots & \ddots & \vdots \\
a_{m1} & a_{m2} & \dots & a_{mn}
\end{pmatrix}.
\]
Each column has sum $k$, since the number of boxes with content $j$ is $k$.
A similar observation gives that the sum of the entries in row $i$ is $(\lambda_i-\mu_i)k$.

The change of variables from GT-patterns to the matrix above is an integral, affine change of variables
and so is the inverse. Thus, lattice points in $\polyP_{\lambdavec/\muvec,\onevec}$
are in bijection with contingency matrices with column sums $1$ and row sums given by $\lambda_i-\mu_i$.

The special case when $\lambda_i-\mu_i=1$ and $n=m$ corresponds to the Birkhoff polytope,
which is the convex hull of all $n \times n$ permutation matrices.
Since $\polyP_{\lambdavec/\muvec,\wvec}$ is compressed,
the Birkhoff polytope is therefore also compressed, which was proved
previously in \cite{StanleyDecompositions80}.
There are several unanswered questions about Birkhoff polytopes, 
such as how to compute their volumes, see \cite{Pak2000}.

\section{Open questions}

We conclude this article with some open questions.
Some of these have been posed in an earlier paper, see \cite{King04stretched}, and we add several new.

\begin{question}
Are all coefficient in the Ehrhart polynomial obtained
from $\polyP_{\lambdavec/\muvec,\wvec}$ non-negative?
\end{question}

\begin{question}
The Gelfand--Tsetlin patterns discussed here are associated with Lie algebras of type $A_n$.
There are polytopes similar to GT-polytopes for other types, see \cite{Berenstein1988,ArdilaBS11}.
Are there similar results regarding the integer decomposition property and compactness for other types of GT-polytopes?
\end{question}

\begin{question}
The numbers $|\Gset_{\lambdavec/\muvec,\wvec}|$ are called (skew) Kostka numbers
or Kostka coefficients, see \cite{Macdonald79symmetric}.
They can be seen as special cases of Littlewood-Richardson coefficients. 
The Littlewood-Richardson coefficients can also be interpreted as the 
number of integer lattice points in certain polytopes, 
for example BZ-polytopes or hive polytopes, \cite{Berenstein92,Knutson99thehoneycomb}.

\emph{What about compactness and IDP among these polytopes?}
This question seems to be related to a conjecture posed by De Loera and McAllister in \cite{LoeraM06},
who conjecture that certain polytopes obtained from the hive conditions have unimodular triangulations.
\end{question}

\begin{question}
The GT-polytopes, hive polytopes, $BZ$-polytopes and the Birkhoff polytopes are all 
polytopes that can (after possibly introducing some slack variables) 
be presented on the form $A\xvec = \yvec$ where $A$ is a matrix with entries in $\{-1,0,1\}$.

\emph{
Is it possible to characterize the matrices $A$,
such that if they are integral, they have the integer decomposition property, or are compact?
}

Note that all totally unimodular matrices give rise to IDP polytopes.
\end{question}

\bibliographystyle{amsplain}

\end{document}